\documentclass[12pt,leqno]{article}
\usepackage{amsmath}
\usepackage{amsfonts,amssymb,latexsym,epsfig,amsthm, amscd}
\usepackage{color}

\makeatletter
\numberwithin{equation}{section}
\makeatother

\newcommand{\E}        {{ {\rm I \hskip -2pt E}}}
\renewcommand{\P}      {{ {\rm I \hskip -2pt P}}}

\newcommand{\R}        {{\rm I \hskip -2pt R}}

\newcommand{\Z}        {{\mathbb Z}}

\newcommand{\M}  {{\mathcal M}}

\newcommand\ta  {a}
\newcommand{\MM} {\overline{\mathcal M}}
\newcommand{\eps}  {{\varepsilon}}

\newcommand{\taub} {\bar{\tau}_b}

\newtheorem{theoreme}{{Theorem}}[section]
\newtheorem{theorem}[theoreme]{{Théeorèem}}
\newtheorem{proposition}[theoreme]{{Proposition}}

\newtheorem{lemma}[theoreme]{Lemma}
\newtheorem{corollary}[theoreme]{Corollary}

\newtheorem{remark}[theoreme]{Remark}

\newcommand{\bpf}[1][Proof]{{\noindent {\sc #1 }}}

\newcommand{\epf}{{\hfill $\square$}}
\begin{document}

\title{Moderate Deviations and Extinction of an Epidemic}
\author{\'E. Pardoux} 

\maketitle

\section*{Abstract}\label{intro}
Consider an epidemic model with a constant flux of susceptibles, in
a situation where the corresponding deterministic epidemic model has a unique stable endemic equilibrium. For the associated stochastic model, whose law of large numbers limit is the deterministic model, the disease free equilibrium is an absorbing state, which is  reached soon or later by the process. However, for a large population size, i.e. when the stochastic model is close to its deterministic limit, 
the time needed for the stochastic perturbations to stop the epidemic may be enormous. In this paper, we discuss how the Central Limit Theorem, Moderate and Large Deviations allow us to give estimates of the extinction time of the epidemic, depending upon the size of the population.

\section{Introduction}
We consider epidemic models where there is a constant flux of susceptible individuals, either because the infected individuals
become susceptible immediately after healing, or after some time during which the individual is immune to the illness, or because there is a constant flux of newborn or immigrant susceptibles.

In the above three cases, for certain values of the parameters,
there is an endemic equilibrium, which is a stable equilibrium of the associated deterministic epidemic model. The deterministic model can be considered as the Law of Large Numbers limit (as the size of the population tends to $\infty$) of a stochastic model, where infections, healings, births and deaths happen according to Poisson processes whose rates depend upon the numbers of individuals in each compartment.  

Since the disease free states are absorbing, it follows from an irreducibility property which is clearly valid in our models, that the epidemic will stop soon or later in the more realistic stochastic model. However, the time which the stochastic perturbances will need to stop the epidemic may be enormous when the size $N$ of the population is large.   
The aim of this paper is to describe, based upon  the Central Limit Theorem, Large and Moderate Deviations, the time it takes for the epidemic to stop in the stochastic model.

The law of large numbers and central limit theorems are rather old. They can be found e.g. in chapter 11 of Ethier and Kurtz \cite{EK}. 
There are also presented, in the framework of epidemic models, in Britton and Pardoux \cite{BP}.
The Large Deviations results are close to those presented in Shwartz and Weiss \cite{SW}, \cite{SW05}, although their assumptions are not quite satisfied in our models. Derivations adapted to our setup can be found  in Kratz and Pardoux \cite{KP},  
Pardoux and Samegni--Kepgnou \cite{PS}, and  Britton and Pardoux \cite{BP}. 
The results concerning moderate deviations are new and constitute the core of this paper. Our derivation is essentially based upon an infinite generalization of the G\"artner--Ellis Theorem, Corollary 4.6.14 from Dembo and Zeitouni \cite{DZ}.  Our main results are Theorem \ref{th:MD} and Theorem \ref{th:WF}. We also give expressions for the rate function in our three models of interest, and in case of the simplest model we give an explicit formula for the quasi--potential. We also compare in that case the upper bound of fluctuations given respectively by the central limit theorem, moderate deviations, and large deviations.

The paper is organized as follows. In section \ref{sec:models}, we describe the three deterministic and stochastic models which we have in mind, namely the SIS, SIRS and SIR model with demography. In section \ref{sec:stochmod}, we give the general formulation of the stochastic models, 
and recall the Law of Large Numbers, the Central Limit Theorem and the Large Deviations, and their application to the time of extinction of an epidemic. In section \ref{sec:moddev}, we establish the moderate deviations result  and explain how it can be used to predict the time taken for an epidemic to cease, depending upon the size of the population. Finally an Appendix establishes an estimate of exponential moments of the integral with respect to a compensated Poisson random measure. This estimate is used several times in our proofs. 

In this paper, the same letter $C$ denotes an arbitrary constant, whose value may change from line to line.

\section{The three models}\label{sec:models}
\subsection{The SIS model}
The deterministic SIS model is the following. Let $s(t)$ (resp. $i(t)$) denote the proportion of susceptible (resp. infectious) individuals in the population. Given an infection parameter $\lambda$, and a recovery parameter $\gamma$, the deterministic SIS model reads
\begin{equation*}
\left\{
\begin{aligned}
s'(t)&=-\lambda s(t)i(t)+\gamma i(t),\\
i'(t)&=\lambda s(t)i(t)-\gamma i(t).
\end{aligned}
\right.
\end{equation*}
Since clearly $s(t)+i(t)\equiv1$, the system can be reduced to a one dimensional ODE. If we let $z(t)=i(t)$, we have $s(t)=1-z(t)$,and 
we obtain the ODE
\[ z'(t)=\lambda z(t)(1-z(t))-\gamma z(t)\,.\]
It is easy to verify that this ODE has a so--called ``disease free equilibrium'', which is $z=0$. If $\lambda>\gamma$, this equilibrium is unstable, and there is an endemic stable equilibrium $z^\ast=1-\gamma/\lambda$.

The corresponding stochastic model is as follows. Let $S^N_t$ (resp. $I^N_t$) denote the proportion of susceptible (resp of infectious)
 individuals in a population of total size $N$.
\begin{equation*}
\left\{
\begin{aligned}
S^N_t&=S^N_0-\frac{1}{N}P_{inf}\left(\lambda N \int_0^tS^N_rI^N_rdr\right)+\frac{1}{N} P_{rec}\left(\gamma N\int_0^t I^N_rdr\right),\\
I^N_t&=I^N_0+\frac{1}{N}P_{inf}\left(\lambda N \int_0^tS^N_rI^N_rdr\right)-\frac{1}{N} P_{rec}\left(\gamma N\int_0^t I^N_rdr\right).
\end{aligned}
\right.
\end{equation*}
Here $P_{inf}(t)$ and $P_{rec}(t)$ are two mutually independent standard  (i.e. rate $1$) Poisson processes. Let us give some explanations, first concerning the modeling, then concerning the mathematical formulation. 

Let $\mathcal{S}^N_t$ (resp. $\mathcal{I}^N_t$) denote the number of susceptible (resp. infectious) individuals in the population.
The equations for those quantities are the above equations, multiplied by $N$. The argument of  $P_{inf}(t)$ reads
\[ \lambda\int_0^t\frac{\mathcal{S}^N_r}{N}\mathcal{I}^N_rdr\,.\]
The formulation of such a rate of infections can be explained as follows. Each infectious individual meets other individuals in the population at some rate $\beta$. The encounter results in a new infection with probability $p$ if the partner of the encounter is susceptible, which happens with probability $S^N_t/N$, since we assume that each individual in the population has the same probability of being that partner, and with probability $0$ if the partner is an infectious individual. Letting $\lambda=\beta p$ and summing over the infectious individuals at time $t$ gives the above rate. Concerning recovery, it is assumed that each infectious individual recovers at rate $\gamma$,
independently of the others. 

\subsection{The SIRS model}
In the SIRS model, contrary to the SIS model, an infectious who heals is first immune to the illness, he is ``recovered'', and only after some time does he loose his immunity and turn to susceptible.  The deterministic SIRS model reads
\begin{equation*}
\left\{
\begin{aligned}
s'(t)&=-\lambda s(t)i(t)+\rho  r(t),\\
i'(t)&=\lambda s(t)i(t)-\gamma i(t),\\
r'(t)&=\gamma i(t)-\rho r(t),
\end{aligned}
\right.
\end{equation*}
while the stochastic SIRS model reads
\begin{equation*}
\left\{
\begin{aligned}
S^N_t&=S^N_0-\frac{1}{N}P_{inf}\left(\lambda N \int_0^tS^N_rI^N_rdr\right)+\frac{1}{N} P_{loim}\left(\rho N\int_0^tR^N_rdr\right),\\
I^N_t&=I^N_0+\frac{1}{N}P_{inf}\left(\lambda N \int_0^tS^N_rI^N_rdr\right)-\frac{1}{N} P_{rec}\left(\gamma N\int_0^t I^N_rdr\right)\\
R^N_t&=R^N_0+\frac{1}{N} P_{rec}\left(\gamma N\int_0^t I^N_rdr\right)-\frac{1}{N}P_{loim}\left(\rho N\int_0^tR^N_rdr\right).
\end{aligned}
\right.
\end{equation*}
These two models could be reduced to two--dimensional models for $z(t)=(i(t),s(t))$ (resp. for $Z^N_t=(I^N_t,S^N_t)$).

\subsection{The SIR model with demography}
In this model, recovered individuals remain immune for ever, but there is a flux of susceptibles by births at a given rate multiplied by $N$, while individuals from each of the three compartments die at rate $\mu$. Thus the deterministic model
\begin{equation*}
\left\{
\begin{aligned}
s'(t) &=\mu -\lambda s(t)i(t) - \mu s(t)
\\
i'(t)&= \lambda s(t)i(t) - \gamma i(t)-\mu i(t) \label{det-SIR-demog}
\\
r'(t) &= \gamma i(t)-\mu r(t),
\end{aligned}
\right.
\end{equation*}
whose stochastic variant reads
\begin{equation*}
\left\{
\begin{aligned}
S^N_t&=S^N_0\!-\!\frac{1}{N}P_{inf}\!\!\left(\!\lambda N \int_0^tS^N_rI^N_rdr\!\right)\!+\!\frac{1}{N} P_{birth}\!\!\left(\mu Nt\right)
\!-\!\frac{1}{N}P_{ds}\!\!\left(\!\mu N\int_0^tS^N_rdr\!\right)\!,\\
I^N_t&=I^N_0+\frac{1}{N}P_{inf}\left(\lambda N \int_0^tS^N_rI^N_rdr\right)-\frac{1}{N} P_{rec}\left(\gamma N\int_0^t I^N_rdr\right)\\
&\qquad-\frac{1}{N} P_{di}\left(\mu N\int_0^tI^N_rdr\right),\\
R^N_t&=R^N_0+\frac{1}{N} P_{rec}\left(\gamma N\int_0^t I^N_rdr\right)-\frac{1}{N} P_{dr}\left(\mu N\int_0^tR^N_rdr\right).
\end{aligned}
\right.
\end{equation*}
\begin{remark} One may think that it would be more natural to decide that births happen at rate $\mu$ times the total population. The total population process would be a critical branching process, which would go extinct in finite time a.s., which we do not want. Next it might seem more natural to replace in the infection rate the ratio $S^N_t/N$ by $S^N_t/(S^N_t+I^N_t+R^N_t)$, which is the actual ratio of susceptibles in the population at time $t$. It is easy to show that $S^N_t+I^N_t+R^N_t$ is close to $N$, so we choose the simplest formulation.
\end{remark}
Again, we can reduce these models to  two--dimensional models for $z(t)=(i(t),s(t))$ (resp. for $Z^N_t=(I^N_t,S^N_t)$), by deleting the $r$ (resp. $R^N$) component.

\section{The stochastic model, LLN, CLT and LD}\label{sec:stochmod}
\subsection{The stochastic model}
The three above stochastic models are of the following form.
\begin{equation}\label{SDE}
\begin{split}
Z^N_t&=z_N+\frac{1}{N}\sum_{j=1}^kh_jP_j\left(N\int_0^t\beta_j(Z^N_s)ds\right)\\
&=z_N+\int_0^tb(Z^N_s)ds+\frac{1}{N}\sum_{j=1}^kh_jM_j\left(N\int_0^t\beta_j(Z^N_s)ds\right),
\end{split}
\end{equation}
where $\{P_j(t),\, t\ge0\}_{0\le j\le k}$ are mutually independent standard Poisson processes, $M_j(t)=P_j(t)-t$, and
$b(z)=\sum_{j=1}^k\beta_j(z) h_j$. $Z^N_t$ takes its values in $\R^d$.

In the case of the SIS model, $d=1$, $k=2$, $h_1=1$, $\beta_1(z)=\lambda z(1-z)$, $h_2=-1$ and $\beta_2(z)=\gamma z$. 

In the case of the SIRS model, $d=2$, $k=3$, $h_1=\begin{pmatrix}1\\-1\end{pmatrix}$, $\beta_1(z)=\lambda z_1z_2$, 
$h_2=\begin{pmatrix}-1\\0\end{pmatrix}$, $\beta_2(z)=\gamma z_1$ and $h_3=\begin{pmatrix}0\\1\end{pmatrix}$, $\beta_3(z)=
\rho(1-z_1-z_2)$.

In the case of the SIR model with demography, we can restrict ourselves to $d=2$, while $k=4$, $h_1=\begin{pmatrix}1\\-1\end{pmatrix}$, $\beta_1(z)=\lambda z_1z_2$, $h_2=\begin{pmatrix}-1\\0\end{pmatrix}$, $\beta_2(z)=(\gamma+\mu) z_1$, $h_3=\begin{pmatrix}0\\1\end{pmatrix}$, 
$\beta_3(z)=\mu$, $h_4=\begin{pmatrix}0\\-1\end{pmatrix}$, $\beta_4(z)=\mu z_2$.

While the above expressions has the advantage of being concise, we shall rather use the following equivalent formulation of \eqref{SDE}. Let $\{\M_j,\, 1\le j\le k\}$ be mutually independent Poisson random measures on $\R_+^2$ with mean measure the Lebesgue measure, and let $\MM_j(ds,du)=\M_j(ds,du)-ds\, du$, $1\le j\le k$. We can rewrite \eqref{SDE} in the form
\begin{equation}\label{SDE+}
\begin{split}
Z^N_t&=z_N+\frac{1}{N}\sum_{j=1}^kh_j\int_0^t\int_0^{N\beta_j(Z^N_s)}\M_j(ds,du)\\
&=z_N+\int_0^tb(Z^N_s)ds+\frac{1}{N}\sum_{j=1}^kh_j\int_0^t\int_0^{N\beta_j(Z^N_s)}\MM_j(ds,du),
\end{split}
\end{equation}
The joint law of $\{Z^N,\, N\ge1\}$ is the same law of a sequence of random elements of the Skorohod space $D([0,T];\R^d)$, whether we use \eqref{SDE} or \eqref{SDE+} for its definition.

Let us state the assumptions which we will need in section \ref{sec:moddev} below. Those are more than necessary for the results of the present section to hold, see \cite{BP} for the proofs.
\begin{align*}
(H.1)&\ \beta_j \ \text{ is bounded }, 1\le j\le k;\\
 (H.2)&\  b\in C^1(\R^d;\R^d),\ \text{ and }\nabla b:\R^d\mapsto \R^{d\times d}\ \text{ is bounded and Lipshitz}.
 \end{align*}

\begin{remark}
In practice, in our models, either the process $Z^N_t$ takes its values in a compact subset of $\R^d$ (this is the case for all models with a constant population size), or else we restrict ourselves to such a situation, by stopping the process when the total population exceeds a given large value, see section 4.2.7 in \cite{BP}. 
\end{remark}

Concerning the initial condition, we assume that for some $z\in[0,1]^d$, $z_N=[Nz]/N$, where
$[Nz]\in\Z^d_+$ is the vector whose $i$--th component is the integer part of the real number $Nz^i$.

\subsection{Law of Large Numbers}
 We have a Law of Large Numbers
\begin{theorem}\label{LLN}
Let $Z^N_t$ denote the solution of the SDE \eqref{SDE}. Then $Z^N_t\to z_t$ a.s. locally uniformly in $t$, where $\{z_t,\, t\ge0\}$ is the unique solution of the ODE
\[  \frac{dz_t}{dt}=b(t,z_t),\quad z_0=x.\]
\end{theorem}
The main argument in the proof of the above theorem is the fact that, locally uniformly in $t$,
\[\frac{P(Nt)}{N}\to t\quad\text{a.s. as }N\to\infty.\]

\subsection{Central Limit Theorem}
We also have a Central Limit Theorem. Let $U^N_t:=\sqrt{N}(Z^N_t-z(t))$.
\begin{theorem}\label{CLT}
As $N\to\infty$, $\{U^N_t,\, t\ge0\}\Rightarrow\{U_t,\, t\ge0\}$ for the topology of locally uniform convergence, where
$\{U_t,\, t\ge0\}$ is a Gaussian process of the form
\begin{equation*}\label{eq:OU}
U_t=\int_0^t\nabla_x b(s,z_s)U_sds+\sum_{j=1}^k h_j\int_0^t\sqrt{\beta_j(s,z_s)}dB_j(s),\ t\ge0\, ,
\end{equation*}
where $\{(B_1(t),B_2(t),\ldots,B_k(t)),\, t\ge0\}$ are mutually independent standard Brownian motions.
\end{theorem}
\subsection{Large Deviations, and extinction of an epidemic}\label{sec:LD}
We denote by $\mathcal{AC}_{T,d}$ the set of absolutely continuous functions from $[0,T]$ into $\R^d$.
For any $\phi\in\mathcal{AC}_{T,d}$, let $\mathcal{A}_{k}(\phi)$ denote the (possibly empty) set of functions $c\in L^1(0,T;\R^k_+)$ such that $c_j(t)=0$ a.e. on the set $\{t,\, \beta_j(\phi_t)=0\}$ and
\begin{equation*}\label{allowed}
  \frac{d\phi_{t}}{dt}=\sum_{j=1}^{k} c_j(t) h_{j}, \quad\text{t a.e}.
\end{equation*}
We define the rate function
\begin{equation*}
I_{T}(\phi) :=
 \begin{cases}
\inf_{c\in\mathcal{A}_{k}(\phi)}I_{T}(\phi|c),& \text{ if } \phi\in\mathcal{AC}_{T,A}; \\
\infty ,& \text{ otherwise.}
\end{cases}
\end{equation*}
where as usual the infimum over an empty set is $+\infty$, and
\begin{equation*}
  I_{T}(\phi|c)=\int_{0}^{T}\sum_{j=1}^{k}g(c_j(t),\beta_{j}(\phi_{t}))dt
\end{equation*}
with $g(\nu,\omega)=\nu\log(\nu/\omega)-\nu+\omega$.
 We assume in the definition of $g(\nu,\omega)$ that for all $\nu>0$, $\log(\nu/0)=\infty$ and $0\log(0/0)=0\log(0)=0$. 
 The collection $Z^N$ obeys a Large Deviations Principle, in the sense that
 \begin{theorem}\label{LD}
 For any open subset $O\subset D([0,T];\R^d)$, 
\[\liminf_{N\to\infty}\frac{1}{N}\log\P\left(Z^{N,z_N}\in O\right)\ge -I_{T,z}(O).\]
For any closed subset $F\subset D([0,T];\R^d)$, 
\[\limsup_{N\to\infty}\frac{1}{N}\log\P(Z^{N,z_N}\in F)\le - I_{T,z}(F)\,.\]
\end{theorem}
A slight reinforcement of this theorem allows us to conclude a Wentzell--Freidlin type of result. 
In what follows, we assume that the first component of $Z^N_t$ (resp. of $z(t)$) is $I^N_t$ (resp. $i(t)$).
Assume that the deterministic ODE which appears in Theorem \ref{LLN} has a unique stable equilibrium $z^\ast$ whose first component satisfies $z^\ast_1>0$. We define 
\[\overline{V}:=\inf_{T>0}\inf_{\phi\in\mathcal{AC}_{T,d}, \phi(0)=z^\ast, \phi_1(T)=0}I_T(\phi).\]
 Let now
 \[ T^{N,z}_{\text{Ext}}=\inf\{t>0,\, Z^N_1(t)=0,\text{ if }Z^N(0)=z_N\}.\]
 We have the
 \begin{theorem}\label{FW}
Given any $\eta>0$, for any $z$ with $z_1>0$,
\begin{equation*}
\lim_{N\to\infty}\mathbb{P}\big(\exp\{N(\overline{V}-\eta)\}<T^{N,z}_{\text{Ext}}<\exp\{N(\overline{V}+\eta)\}\big)=1.
\end{equation*}
Moreover, for all $\eta>0$ and $N$ large enough,
\begin{equation*}
\exp\{N(\overline{V}-\eta)\}\leq\mathbb{E}(T^{N,z}_{\text{Ext}})\leq\exp\{N(\overline{V}+\eta)\}.
\end{equation*}
\end{theorem}
We refer for the proof of this Theorem to \cite{KP} and \cite{BP}.

 It is important to evaluate the quantity $\overline{V}$. Note that it is the value function of an optimal control problem. In case of the SIS model, which is one dimensional, one can solve this control problem explicitly with the help of Pontryagin's maximum principle, see \cite{PBGM}, and deduce in that case that $\overline{V}=\log\frac{\lambda}{\gamma}-1+\frac{\gamma}{\lambda}$.
 For other models, one can compute numerically a good approximation of the value of $\overline{V}$ for each given value of the parameters.
 \subsection{CLT and extinction of an epidemic}
 The discussion of this subsection, which motivates the moderate deviations approach of this paper, is taken from section 4.1 in \cite{BP}.
Consider the SIR with demography. 
\begin{align*}
{i}'(t)&=\lambda i(t) s(t)-\gamma i(t)-\mu i(t),\\
{s}'(t)&= - \lambda i(t) s(t) +\mu -\mu s(t).
\end{align*}
We assume that $\lambda>\gamma+\mu$, in which case there is a unique stable endemic equilibrium,
namely $z^\ast=(i^\ast,s^\ast)=(\frac{\mu}{\gamma+\mu}-\frac{\mu}{\lambda},\frac{\gamma+\mu}{\lambda})$.
We can study the extinction of an epidemic in the above model using the CLT.
We note that the basic reproduction number $R_0$ and the expected relative time of a life an individual is infected, $\varepsilon$, are given by
\begin{equation*}
R_0=\frac{\lambda}{\gamma + \mu} \hskip2cm \varepsilon =\frac{1/(\gamma + \mu)}{1/\mu}=\frac{\mu}{\gamma + \mu}.
\end{equation*}
The rate of recovery $\gamma$ is much larger than the death rate $\mu$ (52 compared to 1/75 for a one week infectious period and 75 year life length) so we use the approximations $R_0\approx \lambda/\gamma$ and $\varepsilon\approx \mu/\gamma$. Denote again by $I^N_t$ the fraction of the population which is infectious in a population of size $N$.
The law of large numbers tells us that for $N$ and $t$ large, $I^N_t$ is close to $i^\ast$. The central limit theorem tell us that $\sqrt{N}(I^N_t-i^\ast)$
converges to a Gaussian process, whose asymptotic variance can be shown to well approximated by $R_0^{-1}$. This suggests that for large $t$, the number of infectious individuals in the population is approximately Gaussian, with mean $Ni^\ast$ and standard deviation $\sqrt{N/R_0}$. If  $Ni^\ast$ and $\sqrt{N/R_0}$ are of the same order, i.e. $N$ is of the same order as 
$\frac{1}{(i^\ast)^2R_0}$, it is likely that the fluctuations described by the central limit theorem explain that the epidemic might cease in time of order one. This gives a critical population size roughly of the order of
\[ N_c\sim\frac{1}{(i^\ast)^2R_0}=\frac{1}{\epsilon^2(1-R_0^{-1})^2R_0},\]
in fact probably a bit larger than that.
% This $N_c$ is rather large since $i^\ast$ is relatively small. Clearly, even if everybody in the population gets ill at some point, being ill one week in a life of average length $75$ years gives a small fraction of infectious in the population. 

Consider measles prior to vaccination. In that case it is known that $R_0\approx 15$, and  $\varepsilon \approx \frac{1/75}{1/(1/52) + 1/75}\approx 1/3750$ we arrive at $N_c\sim (3750)^2/15$, which is almost $10^6$. So, if the population is at most a million (or perhaps a couple of millions), we expect that the disease will go extinct quickly, whereas the disease will become endemic (for a rather long time) in a significantly larger population. This confirms the empirical observation that measles was continuously endemic in UK whereas it died out quickly in Iceland (and was later reintroduced by infectious people visiting the country).

\section{Moderate deviations}\label{sec:moddev}
If the CLT allows to predict extinction of an endemic disease for population sizes under a given threshold $N_c$, and Large Deviations gives predictions for arbitrarily large population sizes, it is fair to look at Moderate Deviations, which describes ranges of  fluctuations between those of the CLT and those of the LD. 

The assumptions $(H.1)$ and $(H.2)$ are assumed to hold throughout this section. 

\subsection{The set--up and preliminary estimates}
We shall use the general model written in the form \eqref{SDE+}.
%Let rewrite our general model \eqref{SDE} in the following equivalent form
%\begin{equation}\label{SDEs}
%\begin{split}
%Z^N_t&=z_N+\frac{1}{N}\sum_{j=1}^k h_j \int_0^t\int_0^{N\beta_j(Z^N_s)}\M_j(ds,du)\\
%&=z_N+\int_0^tb(Z^N_s)ds+\frac{1}{N}\sum_{j=1}^k h_j\int_0^t\int_0^{N\beta_j(Z^N_s)}\MM_j(ds,du),
%\end{split}
%\end{equation}
%where $\M_j,\ j=1,\ldots,k$ are mutually independent Poisson random measures on $\R^2$ with mean measure the Lebesgue
%measure, and for $1\le j\le k$, $\MM_j(ds,du)=\M_j(ds,du)-ds\,du$. From now on, it will be understood that the first coordinate of the process 
%$Z^N_t$ is the proportion of infected individuals in the population.
%
We assume that the limiting law of large numbers ODE
\[ z(t)=z+\int_0^tb(z(s))ds \]
has a unique stable equilibrium point $z^\ast$ such that $z^\ast_1>0$, called the endemic equilibrium, which is such that, provided $z_1(0)>0$, $z(t)\to z^\ast$ as $t\to\infty$.

For the sake of simplifying many formulas below, we chance our coordinates, and let $z^\ast=0$. The reader should be aware of the fact that there is a price to pay for that translation of the origin. Indeed, since in the original coordinate system, the process $Z^N_t$ was living on the set of vectors whose coordinates are integer multiples of $N^{-1}$ (this is essential for the process to remain in the set where it makes sense, i.e. for proportions to remain between $0$ and $1$), the new origin generically does not belong to the set of point in $\R^d$ which our process $Z^N_t$ may visit. The grid on which $Z^N_t$ lives is translated by the vector $z^\ast-\{z^\ast\}_N$, where here and below $\{z\}_N:=[Nz]/N$, $[Nz]$ denoting the vector whose $i$--th component is the integer part of the $i$--th component of $Nz$. However, this minor complexity will
appear only in the formula for the initial condition of the SDE. Once the SDE starts on the correct grid, the solution remains there.

From now on $0$ will be the endemic equilibrium (of course in the translated coordinate system), while $z^\ast\not=0$ will denote that endemic equilibrium in the original coordinates (we shall need it for the formula of the initial condition of the SDE). 

We want to study the moderate deviations at scale $\alpha$ of $Z^N_t$, where $0<\alpha<1/2$. Note that $\alpha=0$ would correspond to the large deviations, and $\alpha=1/2$ to the central limit theorem.
We shall need below to consider the ODE starting from a point close to $z^\ast=0$, namely we shall consider the function $\{z_N(t),\ 0\le t\le T\}$, solution of the ODE
\[ z_N(t)=N^{-\alpha}z+\int_0^tb(z_N(s))ds ,\]
where $z\in\R^d$ is arbitrary. In fact, we shall be more interested in $\overline{z}_N(t):=N^\alpha z_N(t)$, which solves (below we exploit the fact that $b(0)=0$)
\begin{align*}
\overline{z}_N(t)&=z+N^\alpha\int_0^tb(z_N(s))ds\\
&=z+\int_0^t\int_0^1\nabla b(0)\overline{z}_N(s)ds+\int_0^t\int_0^1\left[\nabla b(\theta z_N(s)) -\nabla b(0)\right]d\theta\, \overline{z}_N(s)ds.
\end{align*}
It is not hard to prove that, under our standing assumption $(H.2)$ that $b$ is of class $C^1$ and $\nabla b$ is bounded, as $N\to\infty$, $\overline{z}_N(t)\to\overline{z}(t)$
uniformly for $0\le t\le T$, where $\overline{z}(t)$ solves the linearized ODE near the endemic equilibrium $0$ :
\begin{equation*}\label{linODE}
\overline{z}(t)=z+\int_0^t\nabla b(0)\overline{z}(s)ds\,.
\end{equation*}

We want to study the moderate deviations of the process $Z^N_t$ solution of 
 the SDE \eqref{SDE} with the initial condition $z_N:=\{z^\ast+N^{-\alpha}z\}_N-z^\ast$. This amounts to study the large deviations of 
 $Z^{N,\alpha}_t:=N^\alpha Z^N_t$ at speed $a_N=N^{2\alpha-1}$. We define
 \begin{align*}
 Y^N_t=\frac{1}{N}\sum_{j=1}^k h_j\int_0^t\int_0^{N\beta_j(Z^N_s)}\MM_j(ds,du),\ \text{and }
 Y^{N,\alpha}_t=N^\alpha Y^N_t.
 \end{align*}
  With these notations, the SDE for $Z^{N,\alpha}_t$ reads
  \begin{align*}
  Z^{N,\alpha}_t&= N^\alpha\left(\{z^\ast+N^{-\alpha}z\}_N-z^\ast\right)+\int_0^tN^\alpha b\left(N^{-\alpha}Z^{N,\alpha}_s\right)ds+Y^{N,\alpha}_t
  \\
  &=N^\alpha\left(\{z^\ast+N^{-\alpha}z\}_N-z^\ast\right)+\int_0^t\nabla b(0)Z^{N,\alpha}_sds+\int_0^tV^{N,\alpha}_s ds+Y^{N,\alpha}_t,\ \text{where}
  \\
  V^{N,\alpha}_s&=N^\alpha b\left(N^{-\alpha}Z^{N,\alpha}_s\right)-\nabla b(0)Z^{N,\alpha}_s=
  \left[\int_0^1\left(\nabla b\left(\theta N^{-\alpha}Z^{N,\alpha}_s\right)-\nabla b(0)\right)d\theta\right] Z^{N,\alpha}_s\,.
  \end{align*}

% We now define the difference between the stochastic trajectory and the deterministic trajectory $z_N(\cdot)$ as
%\begin{equation}\label{defX}
%X^N_t:=Z^N_t-z_N(t),\ t\ge0
%\end{equation}
%We moreover define
%\begin{equation}\label{defY}
%Y^N_t:=\frac{1}{N}\sum_{j=1}^k h_j\int_0^t\int_0^{N\beta_j(Z^N_s)}\MM_j(ds,du).
%\end{equation}
%It is clear that
%\begin{equation}\label{eqX}
%X^N_t=[z_N]-z_N+\int_0^t[b(Z^N_s)-b(z_N(s))]ds+Y^N_t.
%\end{equation}
If we let $K:=\sup_z\|\nabla b(z)\|$, we have
\[ \|Z^{N,\alpha}_t\|\le\|z\|+\sqrt{d}N^{\alpha-1}+K\int_0^t\|Z^{N,\alpha}_s\|ds+\|Y^{N,\alpha}_t\| .\]
This combined with Gronwall's Lemma yields
\begin{equation*}\label{estimX}
\| Z^{N,\alpha}_t\|\le e^{Kt}\left(\|z\|+\sqrt{d}N^{\alpha-1}+\sup_{0\le s\le t}\|Y^{N,\alpha}_s\|\right).
\end{equation*}
From the boundedness and Lipschitz property of $\nabla b$, and the formula for $V^{N,\alpha}$, we deduce that
\[\|V^{N,\alpha}_t\|\lesssim\|Z^{N,\alpha}_t\|,\ \text{ and } \|V^{N,\alpha}_t\|\lesssim N^{-\alpha}\|Z^{N,\alpha}_t\|^2\,.\]
We deduce from the last three inequalities
\begin{equation}\label{estimV}
\left\|V^{N,\alpha}_t\right\|\!\lesssim\! N^{\alpha-1}+\left(\|z\|+\sup_{0\le s\le t}\|Y^{N,\alpha}_s\|\right)\wedge
N^{-\alpha}\left(\|z\|+\sup_{0\le s\le t}\|Y^{N,\alpha}_s\|\right)^2\,.
\end{equation}
We now define \[\widetilde{Y}^{N,\alpha}_t=\int_0^tV^{N,\alpha}_sds+Y^{N,\alpha}_t,\ t\ge0,\]
so that 
\begin{equation}\label{eq:ZNa}
Z^{N,\alpha}_t= N^\alpha\left(\{z^\ast+N^{-\alpha}z\}_N-z^\ast\right)+\int_0^t \nabla b(0) Z^{N,\alpha}_sds+
\widetilde{Y}^{N,\alpha}_t\,.
\end{equation}

We will see below that the large deviations of $Z^{N,\alpha}$ will follow from those of $\widetilde{Y}^{N,\alpha}$ by a variant of the contraction principle. 
We first consider the simpler processes
\begin{equation*}\label{defY-}
\overline{Y}^N_t:=\frac{1}{N}\sum_{j=1}^k h_j\int_0^t\int_0^{N\beta_j(0)}\MM_j(ds,du),\ \text{ and } 
\overline{Y}^{N,\alpha}_t=N^\alpha \overline{Y}^N_t
\end{equation*}
which are similar to $Y^N$ and $Y^{N,\alpha}$, but with $Z^N_s$ replaced by $0$. 

\subsection{The limiting logarithmic moment generating function  of $\overline{Y}^{N,\alpha}$}
We note that writing the integral over $[0,N\beta_j(0)]$ as the sum from $\ell=1$ to $\ell=N$ of integrals over $((\ell-1)\beta_j(0),\ell\beta_j(0)]$, we can rewrite $\overline{Y}^{N,\alpha}$ as follows.
\begin{align*}
\overline{Y}^{N,\alpha}&=\frac{1}{N^{1-\alpha}}\sum_{\ell=1}^NQ_\ell(t),\ \text{ where}\\
Q_\ell(t)&=\sum_{j=1}^k h_j\int_0^t\int_{(\ell-1)\beta_j(0)}^{\ell\beta_j(0)}\MM_j(ds,du).
\end{align*}
The processes $Q_1, Q_2,\ldots,Q_N$ are i.i.d., and their law is that of
\begin{equation}\label{eq:Q}
 Q(t)=Q_1(t)=\sum_{j=1}^k h_j\int_0^t\int_{0}^{\beta_j(0)}\MM_j(ds,du).
 \end{equation}
Now let $\nu=(\nu_1,\ldots,\nu_d)$ be a vector of signed measures on $[0,T]$. 
\begin{lemma}\label{MD-}
As $N\to\infty$, (recall that $a_N=N^{2\alpha-1}$)
\[ a_N\log\E\exp\left\{a_N^{-1}\nu(\overline{Y}^{N,\alpha})\right\}\to\frac{1}{2}\E\left(\nu(Q)^2\right).\]
\end{lemma}
\bpf
We use in an essential way the above decomposition of $\overline{Y}^{N,\alpha}$.
\begin{align*}
a_N\Lambda_N(a_N^{-1}\nu)&=N^{2\alpha-1}\log\E\exp\{N^{1-\alpha}\nu(\overline{Y}^N)\} \\
&=N^{2\alpha}\log\E\exp\{N^{-\alpha} \nu(Q)\}\\
&= N^{2\alpha}\log\E\{1+N^{-\alpha} \nu(Q)+N^{-2\alpha} \nu(Q)^2+N^{-3\alpha}R_3^N\}\\
&=N^{2\alpha}\log\left\{1+\frac{N^{-2\alpha}}{2}\E[\nu(Q)^2]+N^{-3\alpha}\E[R_3^N]\right\}\\
&\to\frac{1}{2}\E[\nu(Q)^2],
\end{align*}
provided 
\begin{equation}\label{estimR}
\sup_{N\ge1}|\E[R_3^N]|<\infty,
\end{equation}
 which we will check below. From this it follows that the argument of the logarithm 
on the before last line is greater than or equal to $1$, at least for $N$ large enough, and the final conclusion follows easily from the fact that for any $x\ge0$, $x-x^2/2\le \log(1+x)\le x$. 
Let us now check \eqref{estimR}.  It follows from an exact Taylor formula that
\[|R_3^N|\le \frac{|\nu(Q)|^3}{6}\exp(N^{-\alpha}|\nu(Q|).\]
But $\nu(Q)$ is an affine combination of mutually independent Poisson random variables, so that \eqref{estimR} follows easily
by an explicit computation. \epf

\subsection{The limiting logarithmic moment generating function of $\widetilde{Y}^{N,\alpha}$}
We want to study the large deviations of $\widetilde{Y}^{N,\alpha}$.
The main step will be to prove that Lemma \ref{MD-} remains valid if we replace $\overline{Y}^{N,\alpha}$ by $\widetilde{Y}^{N,\alpha}$, which will follow from the next Proposition.
\begin{proposition}\label{estim-diff}
For any $C>0$, $\nu=(\nu_1,\ldots,\nu_d)$ a vector of signed measures, as $N\to\infty$,
\begin{equation*}\label{neglig}
a_N\log\E\exp\left[C a_N^{-1} \nu(\widetilde{Y}^{N,\alpha}-\overline{Y}^{N,\alpha})\right] \to 0\,.  
\end{equation*}
\end{proposition}
Before we establish that Proposition, let us first prove that it yields the wished result. 
\begin{proposition}\label{pro:correct}
Given Lemma \ref{MD-}, 
if Proposition \ref{estim-diff} holds true, then for any signed measure $\nu$ on $[0,T]$,  as $N\to\infty$,
\[ a_N\log\E\exp\left\{{a_N^{-1} \nu(\widetilde{Y}^{N,\alpha})}\right\}\ \to \frac{1}{2}\E[\nu(Q)^2]\, . \]
 
\end{proposition}
\bpf
For any $\delta>0$, we deduce from H\"older's inequality
\begin{align*}
a_N\log\E&\exp\{a_N^{-1}\nu(\widetilde{Y}^{N,\alpha})\}\\
&=a_N\log\E\left[\exp\{a_N^{-1}\nu(\overline{Y}^{N,\alpha})\}\exp\{a_N^{-1}\nu(\widetilde{Y}^{N,\alpha}-\overline{Y}^{N,\alpha})\}\right] \\
&\le \frac{a_N}{1+\delta}\log\E\exp\{(1+\delta)a_N^{-1}\nu(\overline{Y}^{N,\alpha})\}\\&\quad+
\frac{a_N\delta}{1+\delta}\log\E\exp\left\{\frac{1+\delta}{\delta a_N}\nu(\widetilde{Y}^{N,\alpha}-\overline{Y}^{N,\alpha})\right\},
\end{align*}
so that, if we combine Lemma \ref{MD-} and Proposition \ref{estim-diff}, we deduce that
\[ \limsup_N a_N\log\E\exp\{\nu(a_N^{-1}\widetilde{Y}^{N,\alpha})\}\le\frac{(1+\delta)}{2}\E[\nu(Q)^2],\]
and letting $\delta\to0$, we conclude that
\begin{equation*}\label{limsup}\limsup_N a_N\log\E\exp\{\nu(a_N^{-1}\widetilde{Y}^{N,\alpha})\}\le\frac{1}{2}\E[\nu(Q)^2].
\end{equation*}
For the inequality in the other direction, we note that, by similar arguments,
\begin{align*}
a_N\log\E\exp\left\{\frac{a_N^{-1}}{1+\delta}\nu(\overline{Y}^{N,\alpha})\right\}\le
\frac{a_N}{1+\delta}\log\E\exp\{a_N^{-1}\nu(\widetilde{Y}^{N,\alpha})\}\\
+\frac{a_N\delta}{1+\delta}\log\E\exp\{(\delta a_N)^{-1}\mu(\widetilde{Y}^{N,\alpha}-\overline{Y}^{N,\alpha})\},
\end{align*}
with $\mu=-\nu$, which implies that 
\[ \liminf_N a_N\log\E\exp\{a_N^{-1}\nu(\widetilde{Y}^{N,\alpha})\}\ge\frac{1}{2(1+\delta)}\E[\nu(Q)^2],\]
hence, letting $\delta\to0$ we conclude that
\begin{equation*}\label{liminf}
\liminf_N a_N\log\E\exp\{a_N^{-1}\nu({Y}^{N,\alpha})\}\ge\frac{1}{2}\E[\nu(Q)^2].
\end{equation*}
\epf

The remaining of this subsection will be devoted to the proof of Proposition \ref{estim-diff}.

We note that Proposition \ref{estim-diff} is a consequence of the following two Propositions.
\begin{proposition}\label{Y-barY}
For any $C>0$, as $N\to\infty$,
\begin{equation}\label{neglig_bar}
a_N\log\E\exp\left[C  a_N^{-1}\nu({Y}^{N,\alpha}-\overline{Y}^{N,\alpha})\right] \to 0\,.
\end{equation}
\end{proposition}

\begin{proposition}\label{VN}
For any $C>0$, as $N\to\infty$,
\begin{equation*}
a_N\log\E\exp\left[C  a_N^{-1}\sup_{0\le t\le T}\left\|\int_0^tV^{N,\alpha}_s ds\right\|\right]\to0\,.
\end{equation*}
\end{proposition}

We start with the

\noindent{\sc Proof of Proposition \ref{Y-barY}}
The exponents in the expressions entering \eqref{neglig_bar} are sums over the indices $1\le i\le d$ and $1\le j\le k$. Using repeatedly Schwartz's inequality,  it is sufficient to prove the results with the sum replaced by each of the summands. Therefore in this proof we do as if $d=1$, 
we fix $1\le j\le k$ and for the sake of simplifying the notations, we drop the index $j$.
We note that
\begin{align*}
a_N^{-1}(Y^{N,\alpha}-\overline{Y}^{N,\alpha})&=N^{-\alpha}\int_0^t\int_{N\beta(0)}^{N[\beta(Z^N_s)\vee \beta(0)]}\MM(ds,du)\\
&\quad-N^{-\alpha}\int_0^t\int_{N\beta(Z^N_s)}^{N[\beta(Z^N_s)\vee \beta(0)]}\MM(ds,du)
\end{align*}
It is not hard to see that one can treat each of the two terms on the right separately, and we treat only the first term, the treatment of the second one being quite similar. We note that there exists a compensated standard Poisson process $M(t)$ on $\R_+$ such that the factor of 
$N^{-\alpha}$ in this first term can be rewritten as
\[ W^N_t:=M\left(N\int_0^t(\beta(Z^N_s)-\beta(0))^+ds\right)\,.\]
We need to estimate $\E\exp[C N^{-\alpha}\nu(W^N)]$. If we decompose the signed measure $\nu$ as the difference of two measures as follows
$\nu=\nu_+-\nu_-$, we again  have two terms, and it suffices to treat one of them, say $\nu_+$. Of course it suffices to treat the case where $\nu_+\not=0$. Since the positive constant $C$ is arbitrary, we can w.l.o.g. assume that $\nu_+$ is a probability measure on $[0,T]$. It is then clear that 
\[  \exp\left[CN^{-\alpha}\int_0^TW^N_t\nu_+(dt)\right]\le\exp\left[CN^{-\alpha}\sup_{0\le t\le T}W^N_t\right]\,.\]

We choose a new parameter $0<\gamma<\alpha$, and we write the expression whose expectation needs to be estimated as a sum of two terms as follows.
\begin{equation}\label{decomp}
\begin{split}
\exp&\left\{CN^{-\alpha}\sup_{0\le t\le T}W^N_t\right\}\\
&=\exp\left\{CN^{-\alpha}\sup_{0\le t\le T}W^N_t\right\}{\bf1}_{\sup_{0\le t\le T}\|Z^N_t\|\le N^{-\gamma}} \\
&\quad+\exp\left\{CN^{-\alpha}\sup_{0\le t\le T}W^N_t\right\}{\bf1}_{\sup_{0\le t\le T}\|Z^N_t\|> N^{-\gamma}}.
\end{split}
\end{equation}
We now estimate the first term on the right hand side of \eqref{decomp}. For that sake, we define the stopping time
\[ \sigma_N=\inf\{0\le t\le T;\, \|Z^N_t\|> N^{-\gamma}\}\]
and note that
\begin{align*}
\exp&\left\{CN^{-\alpha}\sup_{0\le t\le T}M\left(N\int_0^t(\beta(Z^N_s)-\beta(0))^+ds\right)\right\}{\bf1}_{\sup_{0\le t\le T}\|Z^N_t\|\le N^{-\gamma}}\\
&\le\exp\left\{CN^{-\alpha}\sup_{0\le t\le T}M\left(N\int_0^{t\wedge\sigma_N}(\beta(Z^N_s)-\beta(0))^+ds\right)\right\}
\end{align*}
Consequently the expectation of the first term on the right of \eqref{decomp} is bounded from above by
\begin{align*}
\E\exp&\left\{CN^{-\alpha}\sup_{0\le t\le T}M\left(N\int_0^{t\wedge\sigma_N}(\beta(Z^N_s)-\beta(0))^+ds\right)\right\}\\
&\le \E\exp\left\{(e^{2CN^{-\alpha}}-1-2CN^{-\alpha})N\int_0^{T\wedge\sigma_N}(\beta(Z^N_t)-\beta(0))^+dt\right\}\\
&\le \exp\left\{CN^{1-2\alpha-\gamma}\right\}\, ,
\end{align*}
where the first inequality follows from  Proposition \ref{le:estimexp} in the Appendix below, and the second one exploits the Lipschitz property of $\beta$.
Consider now the second term on the right hand side of \eqref{decomp}. 
\begin{align*}
\E&\left(\exp\left\{N^{-\alpha}\sup_{0\le t\le T}M\left(N\int_0^t(\beta(Z^N_s)-\beta(0))^+ds\right)\right\}{\bf1}_{\sup_{0\le t\le T}\|Z^N_t\|> N^{-\gamma}}\right)\\
&\le\!\left(\!\E\exp\left\{2N^{-\alpha}\!\!\sup_{0\le t\le T}\!\!M\!\!\left(\!N\!\!\int_0^t\!(\beta(Z^N_s)-\beta(0))^+ds\!\right)\!\right\}\!\right)^{1/2}\\
&\quad\times\P\!\left(\sup_{0\le t\le T}\|Z^N_t\|> N^{-\gamma}\!\right)^{1/2}\\
&\le\exp\left\{CN^{1-2\alpha}\right\}\P\left(\sup_{0\le t\le T}\left\|Y^N_t\right\|>cN^{-\gamma}\right)^{1/2}\, ,
\end{align*}
for some $c,C>0$, where the second inequality follows from Proposition \ref{le:estimexp} and the boundedness of $\beta$.
Estimating the second factor in the last expression amounts to estimating the two probabilities (with another $c>0$)
\begin{equation}\label{estim-prob}
\begin{split}
\P&\left(\sup_{0\le t\le T}M\left(N\int_0^t\beta(Z^N_s)ds\right)>cN^{1-\gamma}\right)\quad\text{and}\\
\P&\left(\sup_{0\le t\le T}\left(-M\left(N\int_0^t\beta(Z^N_s)ds\right)\right)>cN^{1-\gamma}\right)\,.
\end{split}
\end{equation}
We estimate the first probability. For any $\ta>0$,
\begin{equation}\label{majProb}
\begin{split}
\P&\left(\sup_{0\le s\le t}M\left(N\int_0^s\beta(Z^N_r)dr\right)>cN^{1-\gamma}\right)\\
&=\P\left(\sup_{0\le s\le t}\exp\left\{\ta M\left(N\int_0^s\beta(Z^N_r)dr\right)\right\}>\exp\{\ta cN^{1-\gamma}\}\right)\\
&\le e^{-\ta cN^{1-\gamma}}\E\left(\sup_{0\le s\le t}\exp\left\{\ta M\left(N\int_0^s\beta(Z^N_r)dr\right)\right\}\right)\\
&\lesssim e^{-\ta cN^{1-\gamma}} \left(\E\exp\left\{(e^{2\ta}-1-2\ta)N\int_0^t\beta(Z^N_s)ds\right\}\right)^{1/2}\\
&\le \exp\{-\ta cN^{1-\gamma}+(e^{2\ta}-1-2\ta)NCt\}\\
&\le \exp\left\{-\frac{c^2}{8Ct}N^{1-2\gamma}\right\},
\end{split}
\end{equation}
where the second inequality follows from Proposition \ref{le:estimexp} and
the last inequality by optimizing over $\ta>0$. One can easily convince oneself that a similar result holds for the second line 
of \eqref{estim-prob}, making use of  Proposition \ref{le:estimexp} with a negative $a$. Note also for further use that the same result also holds 
in case $\gamma=0$. In that case, the probability on the second line of \eqref{estim-prob} is zero for large enough $c$, in which case the anounced estimate is of course true.

The expectation of the second term of the right hand side of 
\eqref{decomp} is thus dominated by (with $c_1$ and $c_2$ two positive constants)
\[ \exp\{ c_1 N^{1-2\alpha}-c_2 N^{1-2\gamma}\}\to0,\quad\text{ as }N\to\infty.\]
Finally
\begin{align*}
\E\exp&\left\{N^{-\alpha}\sup_{0\le s\le t}M\left(N\int_0^t(\beta(Z^N_s)-\beta(z_s))^+ds\right)\right\}\\
&\le \exp\left\{CN^{1-2\alpha-\gamma}\right\}+\exp\{ c_1 N^{1-2\alpha}-c_2 N^{1-2\gamma}\}
\end{align*}
It follows readily from the inequality $\log(a+b)\le \log(2)+\log(a\vee b)$ that for $N$ large enough
\begin{align*}
a_N\log\E\exp&\left\{N^{-\alpha}\sup_{0\le s\le t}M\left(N\int_0^t(\beta(Z^N_s)-\beta(0))^+ds\right)\right\}\le a_N\log(2)+ C N^{-\gamma},
\end{align*}
which establishes \eqref{neglig_bar}.   \epf
\medskip

We now turn to the second proof.

\noindent{\sc Proof of Proposition \ref{VN}}
Recalling assumption $(H.1)$, we  now define, with $\overline{\beta}_j:=\sup_{z\in\R^d}\beta_j(z)$,
\[ \xi^{N,j}_t:=\frac{1}{N}\int_0^t\int_0^{N\overline{\beta}_j}\M_j(ds,du)\, 1\le j\le k,\]
the event 
\[ A^N_b:=\left\{\sup_{0\le t\le T}\left\|\overline{Y}^N_t\right\|\le b\right\}\bigcap\bigcap_{j=1}^k\left\{\sup_{0\le t\le T}\xi^{N,j}_t\le (1+b')\overline{\beta}_jT\right\}\,,\]
and the stopping time
\[ \taub:=\inf\left\{t>0,\, \left\|\overline{Y}^N_t\right\|>b\right\}\bigwedge\bigwedge_{j=1}^k\inf\left\{t>0,\, \xi^{N,j}_t>(1+b')\overline{\beta}_jT\right\}\, ,\]
where the constant $b>0$ will be chosen below, and the constant $b'>0$ is arbitrary. From the estimate \eqref{estimV},
\begin{align}
a_N\log\E&\left[\exp\left\{a_N^{-1} C \sup_{0\le t\le T}\left\|\int_0^tV^{N,\alpha}_s ds\right\| \right\}\right]\nonumber\\
&\lesssim N^{\alpha-1}+ a_N\log\E\left[\exp\left\{a_N^{-1}C\left(\|z\|+ N^{\alpha}\sup_{0\le t\le T}\|Y^N_t\|\right){\bf1}_{(A^N_b)^c} \right\}\right]\label{term1}
\\&\quad+
a_N\log\E\left[\exp\left\{a_N^{-1}C\left( N^{-\alpha}\|z\|^2+N^\alpha\sup_{0\le t\le T}\|Y^N_s\|^2\right){\bf1}_{A^N_b} \right\}\right],\label{term2}
\end{align}
We take the limit successively in the two terms of the above right hand side.
\noindent{\sc Step 1 : Estimate of \eqref{term1}}
We have
\begin{align*}
\E&\left[\exp\left\{a_N^{-1}C \left(\|z\|+ N^{\alpha}\sup_{0\le t\le T}\|Y^N_t\|\right){\bf1}_{(A^N_b)^c} \right\}\right]\\
%&\le\E\left[\exp\left\{a_N^{-1}C N^{\alpha}\sup_{0\le t\le T}\|Y^N_s\| \right\}{\bf1}_{(A^N_b)^c}
%+{\bf1}_{A^N_b}\right]\\
&\lesssim \exp\left\{CN^{1-2\alpha}\|z\|\right\}\P\left((A^N_b)^c\right)+\E\left[\exp\left\{ CN^{-\alpha}\sup_{0\le t\le T}\|NY^N_s\|\right\}{\bf1}_{(A^N_b)^c}\right] +1,
\end{align*}
We first note that the arguments used in the proof of \eqref{majProb}, in the particular case $\gamma=0$, yield
\begin{align}\label{ProbEventb}
%\begin{split}
\P(\left(A^N_b)^c\right)&\le\P\left(\sup_{0\le t\le T}\left\|\overline{Y}^N_t\right\|> b\right)+
\sum_{j=1}^k\P\left(\sup_{0\le t\le T}(\xi^{N,j}_t-\bar{\beta}t)>b'\bar{\beta}T\right)\nonumber\\
&\lesssim e^{-CN}\, ,
%\end{split}
\end{align}
for some constant $C>0$.
We next estimate the product
\[ \E\left[\exp\left\{ CN^{-\alpha}\sup_{0\le t\le T}\|NY^N_s\|\right\}\right] \P\left((A^N_b)^c\right).\]
For the same reason as in the previous proof, we need only consider the case $d=k=1$.
It follows from Proposition \ref{le:estimexp} that the first factor satisfies
\[\E\left[\exp\left\{ CN^{-\alpha}\sup_{0\le t\le T}\left|\int_0^t\int_0^{N\beta(Z^N_s)}\MM(ds,du)\right|\right\}\right]\lesssim e^{CN^{1-2\alpha}} .\]  
Finally there exist two positive constants $C_1$ and $C_2$ such that
\begin{align*}
\E\left[\exp\left\{a_N^{-1}C \left(\|z\|+ N^{\alpha}\sup_{0\le t\le T}\|Y^N_t\|\right){\bf1}_{(A^N_b)^c} \right\}\right]
&\lesssim 1+\exp\{C_1N^{1-2\alpha}-C_2N\}\\&\le 2\, ,
\end{align*}
for $N$ large enough. So $a_N\log$ of the above tends to $0$, as $N\to\infty$.

\noindent{\sc Step 2 : Estimate of \eqref{term2}}
We first note that 
\begin{align*}
a_N&\log\E\left[\exp\left\{a_N^{-1}C\left( N^{-\alpha}\|z\|^2+N^\alpha\sup_{0\le t\le T}\|Y^N_s\|^2\right){\bf1}_{A^N_b} \right\}\right]\\
&\le CN^{-\alpha}\|z\|^2+a_N\log\E\left[\exp\left\{a_N^{-1}CN^\alpha\sup_{0\le t\le T}\|Y^N_s\|^2{\bf1}_{A^N_b} \right\}\right]\,.
\end{align*}
The first term on the right tends to $0$ as $N\to\infty$. It remains to take care of the second term.
Since $Y^N_t$ is a martingale, it is clear that the process
\[ \left\{\exp\left(a_N^{-1}\frac{C}{2} N^{\alpha} \|Y^N_t\|^2 \right),\ t\ge0\right\} \]
is a submartingale. Consequently, from Doob's $L^2$ submartingale inequality,
 \begin{align*}
\E\left[\sup_{0\le t\le T}\exp\left\{a_N^{-1}C N^{\alpha}\|Y^N_t\|^2{\bf1}_{A^N_b} \right\}\right]
&\le\E\left[\sup_{0\le t\le T\wedge\taub}\exp\left\{a_N^{-1}C N^{\alpha}\|Y^N_t\|^2 \right\}\right]\\
&\le 4\E\left[\exp\left\{a_N^{-1}C N^{\alpha}\|Y^N_{T\wedge\taub}\|^2 \right\}\right]
\end{align*}

Next
\begin{align}\label{major}
%\begin{split}
\E\left[\exp\left\{C N^{1-\alpha}\|Y^N_{T\wedge\taub}\|^2 \right\}\right] &\le
\sqrt{\E\left[\exp\left\{CN^{1-\alpha}\left(\|Y^N_{T\wedge\taub}\|^2-\|\overline{Y}^N_{T\wedge\taub}\|^2\right) \right\}\right]  }\nonumber
\\&\quad\times\sqrt{\E\left[\exp\left\{CN^{1-\alpha}\|\overline{Y}^N_{T\wedge\taub}\|^2 \right\}\right]}
%\end{split}
\end{align}

Consider first the first factor on the right hand side of \eqref{major}. We deduce from the definition of $\taub$ that 
\begin{align*}
\|Y^N_{T\wedge\taub}\|^2-\|\overline{Y}^N_{T\wedge\taub}\|^2&=\left(Y^N_{T\wedge\taub}+\overline{Y}^N_{T\wedge\taub},Y^N_{T\wedge\taub}-\overline{Y}^N_{T\wedge\taub}\right)\\
 &\le (c_{T}+b+2N^{-1}\sup_j\|h_j\|)\left\|Y^N_{T\wedge\taub}-\overline{Y}^N_{T\wedge\taub}\right\|,
 \end{align*}
 with $c_{T}=\sum_{j=1}^k\|h_j\|(1+b')\overline{\beta}_jT$.
 Consequently the square of the first factor on the right of \eqref{major} is bounded from above by 
\begin{align*}
 \E\left[\exp\left\{CN^{1-\alpha}\left\|Y^N_{T\wedge\taub}-\overline{Y}^N_{T\wedge\taub}\right\|\right\}\right]\le
 \E\left[\exp\left\{CN^{1-\alpha}\left\|Y^N_{T}-\overline{Y}^N_{T}\right\|\right\}\right],
 \end{align*}
 where we have used Doob's optional sampling theorem for submartingales. 
 From the same argument as above,we do as if $d=1$, note that
 \begin{align*} 
 \exp\left\{CN^{1-\alpha}\left|Y^N_{T}-\overline{Y}^N_{T}\right|\right\}&\le
 \exp\left\{CN^{1-\alpha}\left(Y^N_{T}-\overline{Y}^N_{T}\right)\right\}\\&\quad+
 \exp\left\{CN^{1-\alpha}\left(\overline{Y}^N_{T}-Y^N_{T}\right)\right\} 
 \end{align*}
 and exploit Proposition \ref{Y-barY} in order to conclude concerning $a_N\log$ of the first factor on the right of \eqref{major}.
 
We next note that
\[ \left\|\overline{Y}^N_{T\wedge\taub}\right\|^2 \le \left\|\overline{Y}^N_{T}\right\|^2{\bf1}_{\{\|\overline{Y}^N_T\|\le b\}} +(b+N^{-1}\sup_j\|h_j\|)^2{\bf1}_{\{\taub<T\}}.\]
Hence the square of the second term on the right of \eqref{major} satisfies
\begin{equation}\label{ineq2}
\begin{split}
\E&\left[\exp\left\{CN^{1-\alpha}\left\|\overline{Y}^N_{T\wedge\taub}\right\|^2 \right\}\right]
\le \sqrt{\E\left[\exp\left\{CN^{1-\alpha} {\bf1}_{\{\taub<T\}}\right\}\right]}\\
&\qquad\qquad\times\sqrt{\E\left[\exp\left\{CN^{1-\alpha}\left\|\overline{Y}^N_{T}\right\|^2{\bf1}_{\left\{\left\|\overline{Y}^N_T\right\|\le b\right\}}\right\}\right]}
\end{split}
\end{equation}
Consider first the second factor on the right of \eqref{ineq2}. We have
\[ \left\|\overline{Y}^N_{T}\right\|^2\le k\sum_{j=1}^k\frac{\|h_j\|^2}{N^2}\left|\int_0^T\int_0^{\beta_j(0)}\MM_j(ds,du)\right|^2\,.\]
Using the Cauchy--Schwartz inequality several times, it is clear that it is sufficient to do as if we had (dropping the index $j$)
\[ \overline{Y}^N_T=\frac{1}{N}\int_0^T\int_0^{N\beta(0)}\MM(ds,du)=\sqrt{a/N}\, \xi_N, \]
with  $a=\beta(0)\, T$ and
$\xi_N = \frac{\theta_N-aN}{\sqrt{aN}}$, where $\theta_N\sim \text{Poi}(aN)$.
We now choose  $b=a/3$. We have

\begin{align*}
\E\exp&\left\{CN^{-\alpha} |\xi_N|^2{\bf1}_{\{|\xi_N|\le\sqrt{aN}/3\}}\right\}\\&=
\sum_{k=\lceil2aN/3\rceil}^{\lfloor4aN/3\rfloor} \exp\left\{CN^{-\alpha}\frac{(k-aN)^2}{aN}\right\}e^{-aN}\frac{(aN)^k}{k!}\\
&\lesssim\int_{-\sqrt{aN}/3}^{\sqrt{aN}/3}  \exp\left\{CN^{-\alpha}x^2\right\}  e^{-aN} \frac{(aN)^{aN+x\sqrt{aN}}}{(aN+x\sqrt{aN})!}\sqrt{aN}dx\\
&\lesssim \frac{1}{\sqrt{2\pi}} \int_{-\sqrt{aN}/3}^{\sqrt{aN}/3}  \exp\left\{CN^{-\alpha}x^2\right\} 
e^{x\sqrt{aN}} \left(1+\frac{x}{\sqrt{aN}}\right)^{-(aN+x\sqrt{aN})} dx
 \\
&\le  \frac{1}{\sqrt{2\pi}} \int_{-\sqrt{aN}/3}^{\sqrt{aN}/3}  
\exp\left\{CN^{-\alpha}x^2\right\} \\
&\quad\times\exp\left\{x\sqrt{aN}-(aN+x\sqrt{aN})\left[\frac{x}{\sqrt{aN}}-\frac{x^2}{2aN}+\frac{x^3}{2(aN)^{3/2}}{\bf1}_{x<0}\right]\right\}dx
\\
&\le  \frac{1}{\sqrt{2\pi}} \int_{-\sqrt{aN}/3}^{\sqrt{aN}/3}  \exp\left\{CN^{-\alpha}x^2-\frac{x^2}{2}+\frac{x^3}{2\sqrt{aN}}{\bf1}_{x>0}\right\}dx
\\
&\le \frac{1}{\sqrt{2\pi}} \int_{-\sqrt{aN}/3}^{\sqrt{aN}/3}  \exp\left\{CN^{-\alpha}x^2-\frac{x^2}{3}\right\}dx
\end{align*}
We have proved that the second factor on the right of \eqref{ineq2} remains bounded, as $N\to\infty$.
We next consider the first factor on the right  of \eqref{ineq2}. We first note that
\[ \exp\left\{4C'N^{1-\alpha}{\bf1}_{\{\tau_{b}<T\}} \right\}\le 1+\exp\left\{4C'N^{1-\alpha}\right\}{\bf1}_{\{\tau_{b}<T\}} \]
But from \eqref{ProbEventb},  $\P\left(\taub<T\right)\lesssim e^{-CN}$.
%\begin{align*}
%\P\left(\taub<T\right)&\le
%\P\left(\sup_{0\le t\le T} |\overline{Y}^N_t|>b\right)+\P\left(\sup_{0\le t\le T} X^N_t>\overline{\beta}T\right)\\
%&\lesssim e^{-cN}.
%\end{align*}

 It follows that the left hand side of \eqref{ineq2} is bounded from above by a constant times
\[ 1+ \exp\{C_1 N^{1-\alpha}-C_2 N\},\]
where $C_1$ and $C_2$ are two positive constants. This last expression is bounded  by $2$, as soon as $N$ is large enough.
Finally $a_N\log$ of the left-hand side of \eqref{ineq2} tends to $0$, as $N\to\infty$.
\epf

\begin{remark}
We note that the full strength of \eqref{estimV} is necessary for the proof of Proposition \ref{VN}. Indeed, while
$a_N\log\E\exp\{CN^{1-\alpha}\sup_{0\le t\le T}\|Y^N_t\|\}$ certainly does not converge to $0$ as $N\to\infty$, clearly
with high probability $\|Y^N_t\|^2$ is smaller than $\|Y^N_t\|$, but $\E\exp\{CN^{1-\alpha}\|Y^N_t\|^2\}=\infty$.
\end{remark}

\subsection{Large deviations of $\widetilde{Y}^{N,\alpha}$}
%The next step consists in establishing exponential tightness of the laws of $\widetilde{Y}^{N,\alpha}$, in the sense that
%\begin{proposition}\label{pro:tensexp}
%For any $R>0$, there exists a compact set $K_R\subset\subset D([0,T];\R^d)$ such that 
%\[ \limsup_N a_N\log \P(\widetilde{Y}^{N,\alpha}\in (K_R)^c)\le-R\,.\]
%\end{proposition}
%\bpf {\color{red}
%The proof of this Proposition follows essentially the lines of the proof of exponential tightness in section 4.2.4 of \cite{BP}, and we will omit it.}
%\epf

We first  define the Fenchel--Legendre transform of 
\begin{align*}
 \Lambda(\nu)&=\frac{1}{2}\E[\nu(Q)^2]\\
 &=\sum_{j=1}^k\frac{\beta_j(0)}{2}\int_{[0,T]^2}s\wedge t\, <h_j,\nu>(ds)<h_j,\nu>(dt),
 \end{align*}
where $Q$ has been defined by \eqref{eq:Q}, $\nu=(\nu_1,\ldots,\nu_d)$ is a vector of signed measures and $<h_j,\nu>(dt)=\sum_{i=1}^d h_j^i\nu_i(dt)$, $h^i_j$ being the $i$--th coordinate of the vector $h_j$. 
We have exploited the fact that $\nu(Q)$ is the sum over $j$ of zero mean mutually independent random variables.
 For each $\phi\in D([0,T];\R^d)$, we define
\[ \Lambda^\ast(\phi)=\sup_{\nu\in(D([0,T];\R^d))^\ast}\{ \nu(\phi)-\Lambda(\nu)\}\,.\] 
The next step will consist in proving that the sequence of processes $\{\widetilde{Y}^{N,\alpha}\}_{N\ge1}$ satisfies a Large Deviation Principle.
%Given Proposition \ref{pro:correct} and Proposition \ref{pro:tensexp}, we deduce from Corollary 4.6.14 from \cite{DZ} the following.
\begin{theorem}\label{th:MD1}
The sequence $\{\widetilde{Y}^{N,\alpha},\, N\ge1\}$ satisfies the Large Deviation Principle in $D([0,T];\R^d)$ equipped with the supnorm topology, with the convex, good rate function $\Lambda^\ast$ and with speed $a_N$, in the sense that for any Borel subset $\Gamma\subset D([0,T];\R^d)$, 
\begin{align*}
-\inf_{\phi\in\mathring{\Gamma}}\Lambda^\ast(\phi)&\le\liminf_N a_N\log\P(\widetilde{Y}^{N,\alpha}\in\Gamma)\\&\le\limsup_Na_N\log\P(\widetilde{Y}^{N,\alpha}\in\Gamma)
\le-\inf_{\phi\in\overline{\Gamma}}\Lambda^\ast(\phi)\,.
\end{align*}
\end{theorem}
Since there is a difficulty with having a topology on  $D([0,T];\R^d)$ which makes it a topological vector space, and allows for a simple characterization of the class of compact sets, we shall use a small detour for the proof of the above Theorem. Recall that
\[ \widetilde{Y}^{N,\alpha}_t=Y^{N,\alpha}_t+\int_0^t V^{N,\alpha}_s ds , \]
where $Y^{N,\alpha}_t$ is piecewise constant, with jumps of size $h_jN^{\alpha-1}$. Let $Y^{N,\alpha,c}_t$ denote the continuous piecewise linear approximation of $Y^{N,\alpha}_t$, which is defined as follows. Let $0=\tau^N_0<\tau^N_1<\tau^N_2<\cdots$ denote the successive jump times of the process $Y^{N,\alpha}_t$. For $i\ge0$, on the interval $[\tau^N_i,\tau^N_{i+1}]$, 
\[ Y^{N,\alpha,c}_t=\frac{\tau^N_{i+1}-t}{\tau^N_{i+1}-\tau^N_i}Y^{N,\alpha}_{\tau^N_i}
+\frac{t-\tau^N_{i}}{\tau^N_{i+1}-\tau^N_i}Y^{N,\alpha}_{\tau^N_{i+1}}\,.\]
Next we define 
\[ \widetilde{\widetilde{Y}}^{N,\alpha}_t=Y^{N,\alpha,c}_t+\int_0^t V^{N,\alpha}_s ds\,.\]
We note that 
\begin{equation}\label{c-equiv} 
\sup_{0\le t\le T}\left\| \widetilde{Y}^{N,\alpha}_t-\widetilde{\widetilde{Y}}^{N,\alpha}_t\right\|\le \sup_j\|h_j\| N^{\alpha-1},
\end{equation}
hence for any $\delta>0$, for $N$ large enough, 
\[ \P\left(\sup_{0\le t\le T}\left\| \widetilde{Y}^{N,\alpha}_t-\widetilde{\widetilde{Y}}^{N,\alpha}_t\right\|\ge\delta\right)=0\,.\]
This implies clearly 
\begin{lemma}\label{le:expequiv}
The two sequences $\{\widetilde{Y}^{N,\alpha}\}_{N\ge1}$ and $\Big\{\widetilde{\widetilde{Y}}^{N,\alpha}\Big\}_{N\ge1}$ are exponentially equivalent in $D([0,T];\R^d)$, equipped with the supnorm topology, in the sense that for each $\delta>0$,
\[ \limsup_{N\to\infty} a_N\log\P\left(\sup_{0\le t\le T}\left\| \widetilde{Y}^{N,\alpha}_t-\widetilde{\widetilde{Y}}^{N,\alpha}_t\right\|\ge\delta\right)=-\infty.\]
\end{lemma}
We shall prove below the following.
\begin{proposition}\label{pro:exptight}
The sequence  $\Big\{\widetilde{\widetilde{Y}}^{N,\alpha}\Big\}_{N\ge1}$ is exponentially tight in $C_0([0,T];\R^d)$, the space of continuous functions 
from $[0,T]$ into $\R^d$, which start from $0$ at $t=0$, in the sense that for any $R>0$, there exists a compact subset 
$K_R\subset\subset C_0([0,T]:\R^d)$ such that
\[\limsup_{N\to\infty} a_N\log\P\left(\widetilde{\widetilde{Y}}^{N,\alpha}\not\in K_R\right)\le -R.\]
\end{proposition} 
Let us now turn to the proof of the above Theorem.

\noindent{\sc Proof of Theorem \ref{th:MD1}}
From \eqref{c-equiv}, we deduce that
\begin{align*}
a_N\log\E\exp&\left\{Ca_N^{-1}\nu\left(\widetilde{Y}^{N,\alpha}_t-\widetilde{\widetilde{Y}}^{N,\alpha}_t\right)\right\}\\
&\lesssim N^{\alpha-1}\to0,
\end{align*}
as $N\to\infty$. Consequently, again by the argument of Proposition \ref{pro:correct}, we deduce from that same Proposition that
for any signed measure $\nu$ on $[0,T]$,  as $N\to\infty$,
\[ a_N\log\E\exp\left\{{a_N^{-1} \nu\left(\widetilde{\widetilde{Y}}^{N,\alpha}\right)}\right\}\ \to \frac{1}{2}\E[\nu(Q)^2]\, . \]
This, together with Proposition \ref{pro:exptight}, allows us to apply Corollary 4.6.14 from \cite{DZ}, to conclude that the sequence 
$\Big\{\widetilde{\widetilde{Y}}^{N,\alpha}\Big\}_{N\ge1}$ satisfies a LDP in $C_0([0,T];\R^d)$ with the good rate function $\Lambda^\ast$, and speed $a_N$. Since $C_0([0,T];\R^d)$ is closed in $D([0,T];\R^d)$ equipped with the supnorm topology, it follows from Lemma 4.1.5 in \cite{DZ} that the same LDP holds in the
latter space, with the same rate function $\Lambda^\ast$, extended to that space by $\Lambda^\ast(\phi)=+\infty$ for
$\phi\in D([0,T];\R^d)\backslash C_0([0,T];\R^d)$. The result now follows from Lemma \ref{le:expequiv}, in view of Theorem 4.2.13 from \cite{DZ}. \epf 

We now turn to the

\noindent{\sc Proof of Proposition \ref{pro:exptight}}
Clearly it suffices to prove both that 
\begin{equation}\label{VNexptight}
\lim_{R\to\infty}\limsup_{N\to\infty}a_N\log \P\left(\sup_{0\le t\le T}\|V^{N,\alpha}_t\|\ge R\right)=-\infty\,,
\end{equation}
and that the sequence $\{Y^{N,\alpha,c}\}_{N\ge1}$ is exponentially tight in $C_0([0,T];\R^d)$. Let us first establish \eqref{VNexptight}.
It follows from \eqref{estimV} that
\[ \| V^{N,\alpha}_t\|\le C\left(\|z\|+1+\sup_{0\le t\le T}\| Y^{N,\alpha}_t\|\right)\,.\]
Consequently, if $R>2C(\|z\|+1)$, with $R'=(2C)^{-1} R$,
\begin{align*}
a_N\log\P(\sup_{0\le t\le T}\|V^{N,\alpha}_t\|>R)&\le a_N\log\P(\sup_{0\le t\le T}\| Y^{N,\alpha}_t\|>R')\\
&\le a_N\log\left(e^{-a_N^{-1}R'}\E\sup_{0\le t\le T}\exp\{a_N^{-1}\|Y^{N,\alpha}_t\|\}\right)\\
&\le -R'+a_N\log \E\left(\sup_{0\le t\le T}\exp\{a_N^{-1}\|Y^{N,\alpha}_t\|\}\right)\,.
\end{align*}
It follows from Doob's submartingale inequality and a combination of Lemma \ref{MD-} and Proposition \ref{Y-barY} that the $\limsup$ as $N\to\infty$ of the second term of the last right hand side is finite. \eqref{VNexptight} clearly follows. 

It remains to consider $Y^{N,\alpha,c}$. Define the modulus of continuity of an element $x\in C_0([0,T];\R^d)$ as $w_x(\delta)=\sup_{0\le s,t\le T, |s-t|\le\delta}\|x(t)-x(s)\|$. It follows from Ascoli's theorem that for any sequence $\{\delta_\ell, \ell\ge1\}$ of positive numbers, the following is a compact subset of $ C_0([0,T];\R^d)$:
\[ \bigcap_{\ell\ge1}\{x:\, w_x(\delta_\ell)\le\ell^{-1}\} \,.\]
Suppose that for each $\ell\ge1$, $R>0$, we can find $\delta_{R,\ell}>0$ such that for all $N\ge1$, 
\begin{equation}\label{YNc-exptight} 
\P(w_{Y^{N,\alpha,c}}(\delta_{R,\ell})\ge\ell^{-1})\le \exp\{-a_N^{-1}(R+\ell)\}\,.
\end{equation}
From this we deduce that 
\begin{align*}
 \P(\cup_{\ell\ge1}\{w_{Y^{N,\alpha,c}}(\delta_{R,\ell})\ge\ell^{-1}\})\le \sum_{\ell\ge1}e^{-a_N^{-1}(R+\ell)}\le e^{-a_N^{-1} R},
 \end{align*}
 so that
 \[ \limsup_{N\to\infty}a_N\log \P(\cup_{\ell\ge1}\{w_{Y^{N,\alpha,c}}(\delta_{R,\ell})\ge\ell^{-1}\})\le -R,\]
 from which the result follows. A sufficient condition for \eqref{YNc-exptight} to be true is that for any $b>0$,
\[ \lim_{\delta\to0}\limsup_{N\to\infty} a_N\log\P\left(w_{Y^{N,\alpha,c}}(\delta)>b\right)=-\infty\,.\]
In turn a sufficient condition for this is that
\begin{equation}\label{YNexp-tight}
\lim_{\delta\to0}\limsup_{N\to\infty} a_N\log\P\left(w_{Y^{N,\alpha}}(\delta)>b\right)=-\infty\,,
\end{equation}
which we now prove. It is not hard to see that
\begin{align*}
\P&\left(w_{Y^{N,\alpha}}(\delta)>b\right)\le2\left(\frac{T}{\delta}+1\right)\sup_{0\le t\le T}\P\left(\sup_{t\le s\le t+2\delta}
\|Y^{N,\alpha}_s-Y^{N,\alpha}_t\|\ge b/2\right)\\
&\le2\left(\frac{T}{\delta}+1\right)\sup_{0\le t\le T}\P\left(\sup_{t\le s\le t+2\delta}\exp\{a_N^{-1}\delta^{-1/2}
\|Y^{N,\alpha}_s-Y^{N,\alpha}_t\|\}\ge \exp\{ba_N^{-1}/2\sqrt{\delta}\}\right)\\
&\le 2\left(\frac{T}{\delta}+1\right)\exp\{-ba_N^{-1}/2\sqrt{\delta}\}\sup_{0\le t\le T}\E\exp\{a_N^{-1}\delta^{-1/2}\|Y^{N,\alpha}_{t+2\delta}-Y^{N,\alpha}_t\|\},
\end{align*}
where we have used Doob's submartingale inequality at the last step. Clearly
\begin{align*}
\exp\{a_N^{-1}\delta^{-1/2}\|Y^{N,\alpha}_{t+2\delta}-Y^{N,\alpha}_t\|\}\le\prod_{j=1}^k\exp\left\{N^{-\alpha}\delta^{-1/2}\| h_j\|\left|\int_t^{t+2\delta}\int_0^{N\beta_j(Z^N_s)}\MM_j(ds,du)\right|\right\}
\end{align*}
Using repeatedly Cauchy--Schwartz's inequality, we see that it suffices to estimate for each $j$
\begin{align*}
\E\exp&
\left\{CN^{-\alpha}\delta^{-1/2}\left|\int_t^{t+2\delta}\int_0^{N\beta_j(Z^N_s)}\MM_j(ds,du)\right|\right\}
\\
&\le 
\E\exp\left\{CN^{-\alpha}\delta^{-1/2}\int_t^{t+2\delta}\int_0^{N\beta_j(Z^N_s)}\MM_j(ds,du)\right\}
\\
&+ \E\exp\left\{-CN^{-\alpha}\delta^{-1/2}\int_t^{t+2\delta}\int_0^{N\beta_j(Z^N_s)}\MM_j(ds,du)\right\}\\
&\le 2\exp\{8 C^2N^{1-2\alpha}\bar{\beta_j}\},
\end{align*}
where $\bar{\beta}_j=\sup_{z}\beta_j(z)$, we have used Proposition \ref{le:estimexp} and the inequality $e^x-1-x\le x^2$, valid for $x\le\log(2)$, which we have applied with $x=2CN^{-\alpha}\delta^{-1/2}$ and $x=-2CN^{-\alpha}\delta^{-1/2}$ (recall that we will first let $N\to\infty$). Putting together the last estimates yields
\[ \limsup_{N\to\infty}a_N\P\left(w_{Y^{N,\alpha}}(\delta)>b\right)\le -\frac{b}{2\sqrt{\delta}} + C.\]
\eqref{YNexp-tight} follows, and the Proposition is proved. \epf

\subsection{Computation of the rate function $\Lambda^\ast$}
Let us compute $\Lambda^\ast$ in the three examples which we discussed above in section \ref{sec:models}. 
Here we do not translate $z^\ast$ to the origin.

\subsubsection{Computation of $\Lambda^\ast$ for the SIS model}
Recall that in this case $d=1$, $k=2$, $h_1=1$, $\beta_1(z)=\lambda z(1-z)$, $h_2=-1$, $\beta_2(z)=\gamma z$.
If $\lambda>\gamma$, there is a unique stable endemic equilibrium $z^\ast=1-\gamma/\lambda$.
We first compute 
\[\Lambda(\nu)=\frac{1}{2}\E\int_{[0,T]\times[0,T]}Q(s)Q(t)\nu(ds)\nu(dt),\]
where 
\[ Q(t)=\int_0^t\int_0^{\beta_1(z^\ast)}\MM_1(ds,du)-\int_0^t\int_0^{\beta_2(z^\ast)}\MM_2(ds,du)\,.\]
It is easy to check that $\E[Q(t)Q(s)]=\sigma^2(z^\ast)\,{s\wedge t}$, where 
\[ \sigma^2(z^\ast)=\beta_1(z^\ast)+\beta_2(z^\ast)=2\frac{\gamma}{\lambda}(\lambda-\gamma).\]
 Consequently 
\[ \Lambda(\nu)=\frac{\sigma^2(z^\ast)}{2}\int_{[0,T]\times[0,T]}{s\wedge t}\ \nu(ds)\nu(dt)\,.\]

We now need to compute $\Lambda^\ast(\phi)$ in case $\phi\in C^2([0,T])$. We should take the supremum over the signed measures $\nu$ on $[0,T]$ of the quantity
\[ \int_{[0,T]}\phi(t)\nu(dt)-\frac{\sigma^2(z^\ast)}{2} \int_{[0,T]\times[0,T]}{s\wedge t}\, \nu(ds)\nu(dt)\,.\]
The supremum is achieved at the signed measure $\nu$ which makes the gradient with respect to $\nu$ of the above zero, if any.
We first note that for such a $\nu$ to exist, we need that $\phi(0)=0$, unless $\Lambda^\ast(\phi)=+\infty$. Now the optimal $\nu$ must satisfy
\begin{align*}
 \phi(t)&=\sigma^2(z^\ast)\int_{[0,T]}{s\wedge t}\,\nu(ds)\\
 &=\sigma^2(z^\ast)\int_{[0,t]} s\,\nu(ds)+ \sigma^2(z^\ast)t \int_{(t,T]}\nu(dt).
 \end{align*}
So necessarily
 \[ \nu(dt)=-\frac{\phi''(t)}{\sigma^2(z^\ast)}dt+\frac{\phi'(T)}{\sigma^2(z^\ast)}\delta_T(dt).\]
 Substituting this signed measure $\nu$ in the above formula, we obtain that
  \begin{align*}
   \int_{[0,T]}\phi(t)\nu(dt)&=\frac{\phi'(T)}{\sigma^2(z^\ast)}\phi(T)-\int_0^T\frac{\phi''}{\sigma^2(z^\ast)}(t)\phi(t)dt\\
   &=\frac{1}{\sigma^2(z^\ast)}\int_0^T|\phi'(t)|^2dt\, .
   \end{align*}
   Consequently
 \begin{align*}
 \Lambda(\phi)=\begin{cases}
 \frac{1}{2\sigma^2(z^\ast)}\int_0^T|\phi'(t)|^2dt,&\text{if $\phi(0)=0$ and $\phi$ is absolutely continuous;}\\
    +\infty,&\text{otherwise}.
    \end{cases}
    \end{align*}
 
\subsubsection{Computation of $\Lambda^\ast$ for the SIRS model}
In this model, $d=2$ and $k=3$. We have $h_1=\binom{1}{-1}$, $\beta_1(z)=\lambda z_1z_2$, $h_2=\binom{-1}{0}$, $\beta_2(z)=\gamma z_1$, and $h_3=\binom{0}{1}$, $\beta_3(z)= \rho(1-z_1-z_2)$. In the case $\lambda>\gamma$, there is a unique stable endemic equilibrium, namely
$z^\ast=\binom{\frac{\rho}{\gamma+\rho}\left(1-\frac{\gamma}{\lambda}\right)}{\frac{\gamma}{\lambda}}$. In order to simplify the notations, we shall write $a=\beta_1(z^\ast)$, $b=\beta_2(z^\ast)$, $c=\beta_3(z^\ast)$  and $A=ab+ac+bc$. We  have
\begin{align*}
\Lambda(\nu)&=\frac{a}{2}\int_{[0,T]\times[0,T]}s\wedge t(\nu_1-\nu_2)(ds)(\nu_1-\nu_2)(dt)\\
&\quad+\frac{b}{2}\int_{[0,T]\times[0,T]}s\wedge t\nu_1(ds)\nu_1(dt)+\frac{c}{2}\int_{[0,T]\times[0,T]}s\wedge t\nu_2(ds)\nu_2(dt)
\end{align*}
The functional to be maximized with respect to $\nu$ if
\begin{align*} 
<\nu_1,\phi_1>&+<\nu_2,\phi_2>-\frac{a}{2}\int_{[0,T]\times[0,T]}s\wedge t(\nu_1-\nu_2)(ds)(\nu_1-\nu_2)(dt)\\&
-\frac{b}{2}\int_{[0,T]\times[0,T]}s\wedge t\nu_1(ds)\nu_1(dt)-\frac{c}{2}\int_{[0,T]\times[0,T]}s\wedge t\nu_2(ds)\nu_2(dt)
\end{align*}
Writing that the gradient w.r.t. $\nu_1$ and $\nu_2$ of this functional is zero leads to the identities
\begin{align*}
\phi_1(t)&=a\int_{[0,T]}s\wedge t (\nu_1-\nu_2)(ds)+b\int_{[0,T]}s\wedge t \nu_1(ds),\\
\phi_2(t)&=a\int_{[0,T]}s\wedge t (\nu_2-\nu_1)(ds)+c\int_{[0,T]}s\wedge t \nu_2(ds).
\end{align*}
This implies the identities
\begin{align*}
\nu_1(dt)&=-\left(\frac{a+c}{A}\phi''_1(t)+\frac{a}{A}\phi''_2(t)\right)dt+\left(\frac{a+c}{A}\phi'_1(T)+\frac{a}{A}\phi'_2(T)\right)\delta_T,\\
\nu_2(dt)&=-\left(\frac{a}{A}\phi''_1(t)+\frac{a+b}{A}\phi''_2(t)\right)dt+\left(\frac{a}{A}\phi'_1(T)+\frac{a+b}{A}\phi'_2(T)\right)\delta_T.
\end{align*}
Finally we deduce that $\Lambda^\ast(\phi)$ is $+\infty$ unless $\phi$ is absolutely continuous and $\phi(0)=0$, in which case
\[†\Lambda^\ast(\phi)=\frac{a}{2A}\int_0^T|\phi'_1(t)+\phi'_2(t)|^2dt
+\frac{c}{2A}\int_0^T|\phi'_1(t)|^2dt
+\frac{b}{2A}\int_0^T|\phi'_2(t)|^2dt\, . \]

\subsubsection{Computation of $\Lambda^\ast$ for the SIR model with demography}
In this case, $d=2$, $k=4$, $h_1=\binom{1}{-1}$, $\beta_1(z)=\lambda z_1z_2$, $h_2=\binom{-1}{0}$, $\beta_2(z)=(\gamma+\mu)z_1$,
$h_3=\binom{0}{1}$, $\beta_3(z)=\mu$ and $h_4=\binom{0}{-1}$, $\beta_4(z)=\mu z_2$.
In the case $\lambda>\gamma+\mu$, there is a unique stable endemic equilibrium, namely $z^\ast=\binom{\mu\left(\frac{1}{\gamma+\mu}-\frac{1}{\lambda}\right)}{\frac{\gamma+\mu}{\lambda}}$. We shall use the notations
$a=\beta_1(z^\ast)$, $b=\beta_2(z^\ast)$, $c=\beta_3(z^\ast)+\beta_4(z^\ast)$  and $A=ab+ac+bc$
We have
\begin{align*}
\Lambda(\nu)&=\frac{a}{2}\int_{[0,T]\times[0,T]}s\wedge t(\nu_1-\nu_2)(ds)(\nu_1-\nu_2)(dt)\\
&\quad+\frac{b}{2}\int_{[0,T]\times[0,T]}s\wedge t\nu_1(ds)\nu_1(dt)+\frac{c}{2}\int_{[0,T]\times[0,T]}s\wedge t\nu_2(ds)\nu_2(dt)
\end{align*}
Formally the functional $\Lambda(\nu)$ has exactly the same form as in the case of the SIRS model, only the constants have different values. The same computations as in the previous subsection lead to the same result, namely that
$\Lambda^\ast(\phi)$ is $+\infty$ unless $\phi$ is absolutely continuous and $\phi(0)=0$, in which case
\[\Lambda^\ast(\phi)=\frac{a}{2A}\int_0^T|\phi'_1(t)+\phi'_2(t)|^2dt
+\frac{c}{2A}\int_0^T|\phi'_1(t)|^2dt
+\frac{b}{2A}\int_0^T|\phi'_2(t)|^2dt\, . \]

\subsection{Moderate deviations of $Z^N$}
We again equip $D([0,T];\R^d)$ with the supnorm topology.
Let for $z\in\R^d$ $F_z:D([0,T];\R^d)\mapsto D([0,T];\R^d)$ be the continuous map which to $x$ associates $y$ solution of the ODE
\[ y(t)=z+\int_0^t \nabla b(0)y(s)ds +x(t),\]
and for each $N\ge1$ $F_{z,N}:D([0,T];\R^d)\mapsto D([0,T];\R^d)$ be the continuous map which to $x$ associates $y_N$
solution of the ODE
\[ y_N(t)=N^\alpha(\{z^\ast+N^{-\alpha}z\}_N-z^\ast)+\int_0^t \nabla b(0) y_N(s)ds +x(t).\]
We have
\begin{equation}\label{exp-tight}
 F_{z,N}(x)(t)-F_z(x)(t)=\exp[\nabla b(0) t]\, \left(N^\alpha(\{z^\ast+N^{-\alpha}z\}_N-z^\ast)-z\right),
 \end{equation}
which converges to $0$ as $N\to\infty$, uniformly in $t\in[0,T]$ and $x\in D([0,T];\R^d)$. 
%where $A(t)$is the $d\times d$ matrix valued function of $t$ which solves the ODE
%\[ \frac{dA}{dt}(t)=\nabla b(z_t) A(t),\ A(0)=I_d,\]
%where $I_d$ denotes here the $d\times d$ identity matrix.
%Hence there exists $C>0$ such that $\|F(x)-F_N(x)\|_T\le C N^{\alpha-1}$. 
We want to study the moderate deviations of $Z^N$, or in other words the large deviations of $Z^{N,\alpha}=N^\alpha Z^N$. 
In what follows, we shall denote by 
$Z^{N,\alpha}_z$ the process $Z^{N,\alpha}$ starting from $Z^{N,\alpha}(0)=N^\alpha(\{z^\ast+N^{-\alpha}z\}_N-z^\ast)$. 
From \eqref{eq:ZNa}, $Z^{N,\alpha}_z=F_{z,N}(\widetilde{Y}^{N,\alpha})$, hence
the following statement is a consequence of Theorem \ref{th:MD1}, \eqref{exp-tight} and Corollary 4.2.21 from \cite{DZ}.

\begin{theorem}\label{th:MD} Assume that $(H.1)$ and $(H.2)$ hold.
The collection of processes $\{Z^{N,\alpha}_z(t),\ 0\le t\le T\}_{N\ge1}$ satisfies a large deviations principle with speed $a_N=N^{2\alpha-1}$ and the good rate function
\begin{align*} 
I_{z,T}(\phi)&=\Lambda^\ast(F_z^{-1}(\phi))\\
&=\begin{cases} \Lambda^\ast\left(\phi(\cdot)-z-\nabla b(0)\int_0^\cdot\phi(s)ds\right)&\text{if $\phi(0)=z$};\\
+\infty,&\text{otherwise.}\end{cases} 
\end{align*}
More precisely, for any Borel subset $\Gamma\subset D([0,T];\R^d)$,
\begin{align*}
-\inf_{\phi\in\mathring{\Gamma}}I_{z,T}(\phi)&\le\liminf_N a_N\log\P(Z^{N,\alpha}_{z} \in\Gamma)\\&\le
\limsup_Na_N\log\P(Z^{N,\alpha}_{z} \in\Gamma)
\le-\inf_{\phi\in\overline{\Gamma}} I_{z,T}(\phi)\,.
\end{align*}
\end{theorem}

Since the mapping $F_z$ has the nice property that $F_z(x)(t)-F_{z'}(x)(t)=\exp[\nabla b(0) t](x-x')$, it follows readily again from 
Corollary 4.2.21 in \cite{DZ} that the above result can be extended to the following statement. 
\begin{theorem}\label{th:MDext}
Assume that $(H.1)$ and $(H.2)$ hold.
For any closed set $F\subset D([0,T];\R^d)$, for any sequence $z_N\to z$, 
\[ \limsup_{N\to\infty}a_N\log\P(Z^{N,\alpha}_{z_N}\in F)\le-\inf_{\phi\in F}I_{z,T}(\phi)\,. \]
For any open set $G\subset D([0,T];\R^d)$,  for any sequence $z_N\to z$,
\[ \liminf_{N\to\infty}a_N\log\P(Z^{N,\alpha}_{z_N}\in G)\ge-\inf_{\phi\in G}I_{z,T}(\phi)\,.\]
\end{theorem}
From this last Theorem, we can deduce, with the same proof as that of Corollary 5.6.15 in \cite{DZ}, the following Corollary.
\begin{corollary}\label{cor:LD}
Assume that $(H.1)$ and $(H.2)$ hold.
Let $K$ denote an arbitrary compact subset of $\R^d$.

For any closed set $F\subset D([0,T];\R^d)$,  
\[ \limsup_{N\to\infty}a_N\log\sup_{z\in K}\P(Z^{N,\alpha}_{z}\in F)\le-\inf_{\phi\in F, z\in K}I_{z,T}(\phi)\,. \]
For any open set $G\subset D([0,T];\R^d)$, 
\[ \liminf_{N\to\infty}a_N\log\inf_{z\in K}\P(Z^{N,\alpha}_{z}\in G)\ge-\sup_{z\in K}\inf_{\phi\in G}I_{z,T}(\phi)\,.\]

\end{corollary}

\subsection{Wentzell--Freidlin theory and extinction of an epidemic}

We now define
\begin{align*}
V(z,z',t)&=\inf_{\phi, \, \phi(0)=z,\phi(t)=z'}I_{z,t}(\phi),\\
V(z,z')&=\inf_{t>0}V(z,z',t),\\
\overline{V}_a&=\inf_{z,\, z_1= -a}V(0,z),
\end{align*}
where $a>0$, and we recall that we have translated the endemic equilibrium $z^\ast$ at the origin.

We can now state our main result.
\begin{theorem}\label{th:WF}
Assume that $(H.1)$ and $(H.2)$ hold. For some $a>0$, let $T^N_{z,a}:=\inf\{t>0,\ Z^{N,\alpha}_{z,1}(t)\le -a\}$, where $Z^{N,\alpha}_{z,1}(t)$ denotes the first coordinate of the process $Z^{N,\alpha}_{z}(t)$. The following hold.

For any $z\in\R^d$ such that $z_1>-a$, and any $\eta>0$,
\[\lim_{N\to\infty}\P\left(e^{a_N^{-1}(\overline{V}_a-\eta)} < T^N_{z,a} < e^{a_N^{-1}(\overline{V}_a+\eta)}\right)=1\,,\]
and
\[\lim_{N\to\infty} a_N\log \E(T^N_{z,a})=\overline{V}_a\,.\]
\end{theorem}

%
%We want to conclude from the Wentzell--Freidlin theory an estimate of the time needed for the first coordinate of $N^\alpha(Z^N_t-z^\ast)$ to make a deviation of $-c$, i.e. to go from $0$ to $-c$, which, for the value of $N$ such that $z_1^\ast=cN^{-\alpha}$, means the time for the first coordinate of $Z^N_t$ to hit $0$ (in the non translated coordinate system). For that sake, we first compute
%\[\overline{V}_c=\min_{T>0}\min_{\phi,\ \phi(0)=0, \phi(T)=-c}I_T(\phi).\]
%An application of Pontryaguin's maximum principle, see \cite{PBGM}, yields
%\[ \overline{V}_c=\frac{\lambda-\gamma}{\sigma^2} c^2\, .\]
%$a=\beta_1(z^\ast)$, $b=\beta_2(z^\ast)$, $c=\beta_3(z^\ast)$, $d=\beta_4(z^\ast)$. 
%Using the same arguments as in \cite{KP} and \cite{BP}, we then deduce from Theorem \ref{th:MD1}
%\begin{theorem}
%Let $T^{N,\alpha}_c:=\inf\{t>0,\, \widetilde{Z}^{N,\alpha}_t\le -c\}$. For any $\delta>0$, 
%\[ \lim_{N\to\infty}\P\left(\exp\{a_N^{-1}(\overline{V}_c-\delta)\}<T^{N,\alpha}_c<\exp\{a_N^{-1}(\overline{V}_c+\delta)\}\right)=1\, .\]
%Moreover
%\[ \lim_{N\to\infty}a_N\E(T^{N,\alpha}_c)=\overline{V}_c\,.\]
%\end{theorem}

Given Corollary \ref{cor:LD}, the proof of the above result follows the exact same steps as that of Theorem 5.7.11 in \cite{DZ}, with some minor modifications, to adapt to the fact that our processes have discontinuous trajectories, see the proof of Theorem 7.14 in \cite{KP}, or of Theorem 4.2.17 in \cite{BP}.

Recall that $a_N^{-1}=N^{1-2\alpha}$. In the CLT regime, $\alpha=1/2$, $a_N^{-1}=1$, while in the LD regime, $\alpha=0$, 
$a_N^{-1}=N$.
\subsubsection{Interpretation. The critical population size}
Going back to the original coordinates, i.e. $z^\ast\not=0$, we should interpret $Z^{N,\alpha}(t)$ as $Z^{N,\alpha}(t)=N^\alpha(Z^N(t)-z^\ast)$. So (dropping the index for the starting point in order to simplify our notations), $T^N_a$ is the first time  when
$Z^N_1(t)\le z^\ast_1-aN^{-\alpha}$. For $T^N_a$ to be finite, we need to have $z^\ast_1-aN^{-\alpha}\ge0$, since $Z^N_1(t)$ cannot become negative. This is of course no problem for the limit theorem, since $aN^{-\alpha}\to0$ as $N\to\infty$, while $z^\ast_1$ is fixed. 
However, a deviation of the order of $-aN^{-\alpha}$ is enough for $Z^N_1(t)$ to hit zero, if $z^\ast_1$ is of the order of $N^{-\alpha}$, which means that $N$ is of the order of $(z^\ast_1)^{-1/\alpha}$. 
%Let us now compute the corresponding critical population size. 
$e^{N^{1-2\alpha}\overline{V}_{a}}$ is the order of magnitude of the time needed for $Z^N_t-z_t$ to make a deviation of size $aN^{-\alpha}$. This is sufficient to extinguish an epidemic, provided $z_1^\ast$ is of the same order, so that the corresponding critical size is $N_{c,\alpha}\sim (1/z_1^\ast)^{1/\alpha}$, which is roughly the CLT critical population size raised to the power $1/2\alpha$. 

\subsubsection{The value of $\overline{V}_a$ in the SIS model} 
In the particular case of the SIS model, we can compute explicitly the value of the quasi--potential $\overline{V}_a$. In this case, $d=1$, the linearized ODE around the endemic equilibrium translated at $0$ reads
\[ \dot{x}=-(\lambda-\gamma) x+u,\]
and the cost functional to minimize is 
\[ I_T(u)=\frac{\lambda}{4\gamma(\lambda-\gamma)}\int_0^Tu(t)^2dt\,.\]
We are looking for the minimal cost for driving $x$ from $0$ to $-a$. We now exploit the Pontryagin maximum principle, see \cite{PBGM}.
The Hamiltonian reads
\[ H(x,p,u)=-(\lambda-\gamma)px+pu-\frac{\lambda}{4\gamma(\lambda-\gamma)}u^2\,.\]
The optimal control $\hat{u}$ must maximize the Hamiltonian, so it satisfies $\hat{u}=\frac{4\gamma(\lambda-\gamma)}{\lambda}p$. Since the final time is free and the system is autonomous, the Hamiltonian vanishes along the optimal trajectory, so that along such a trajectory, either $p=0$, in which case $\hat{u}=0$, or else  $x=\frac{\gamma}{\lambda}p$, hence  
$\hat{u}=2(\lambda-\gamma)x$. Finally the pieces of optimal trajectory which move towards the origin correspond to $u\equiv0$, those which move away from the origin (this is the case we are interested in) satisfy the time reversed ODE $\dot{x}=(\lambda-\gamma)x$. There is no optimal trajectory from $x=0$ to $x=-a$. However, if we start from $x=-\eps$, the optimal trajectory is $x(t)=-e^{(\lambda-\gamma)t}\eps$, so $\hat{u}(t)=-2(\lambda-\gamma)e^{(\lambda-\gamma)t}\eps$, the final state $-a$ is reached at time 
$(\lambda-\gamma)^{-1}\log(a/\eps)$, and the optimal cost is $\frac{\lambda}{2\gamma}(a^2-\eps^2)$. 
A possible sub--optimal control starting from $0$ is as follows. Choose $u=-1$ for a time of order $\eps$, until $x(t)$ reaches $-\eps$, 
whose cost is of the order of $\eps$, and then choose the optimal feedback, until $-a$ is reached.
Letting $\eps\to0$, the total cost converges to
\[ \overline{V}_a=\frac{\lambda}{\gamma}\frac{a^2}{2}\,.\]

\subsection{Comparison between the  CLT, MD and LD} We do that comparison in case of the SIS model, for which we have explicit expressions for the rate functions and the quasi--potentials. We still translate $z^\ast$ at the origin, and start our process at the origin : $Z^N_0=0$. To make a change with the above, we fix $a>0$ and want to compare (for $t$ large) the upper bounds for $\P(N^\alpha Z^N_t\ge a)$ in the three cases $\alpha=1/2$ (the central limit theorem), $0<\alpha<1/2$ (moderate deviations) and $\alpha=0$ (large deviations).

We start with the {\bf central limit theorem}.
It is easy to see that $U_t=\lim_{N\to\infty} \sqrt{N}Z^N_t$ solves the SDE
\[ U_t=-(\lambda-\gamma)\int_0^t U_s ds + \sqrt{2\gamma(\lambda-\gamma)/\lambda} B_t,\]
so that the asymptotic variance of $U_t$ is $\gamma/\lambda$. Consequently for $a>0$ fixed and any $\eta>0$, there exist $t$ and $N$ large enough such that we have the following upper bound for the probability of a positive deviation of $\sqrt{N}Z^N_t$ 
\[ \P(\sqrt{N}Z^N_t\ge a)\le \exp\left\{-\frac{\lambda a^2}{2\gamma}+\eta\right\}\,.\]
This bound follows from the following estimate, valid for a $\mathcal{N}(0,\sigma^2)$ random variable $\xi$: 
$\P(\xi>a)\le e^{-\lambda a}e^{\frac{\lambda^2\sigma^2}{2}}$ after optimizing over $\lambda>0$.

Consider next the {\bf moderate deviations}.
Theorem \ref{th:MD} combined with the computation from the last subsection indicates that for $0<\alpha<1/2$, any $\eta>0$, there exists $t$ and $N$ large enough such that 
\[ \P(N^\alpha Z^N_t\ge a)\le \exp\left\{-N^{1-2\alpha}\left(\frac{\lambda a^2}{2\gamma}-\eta\right)\right\}\,.\]

We finally consider  the {\bf large deviations}. Here we need to assume that $a< \gamma/\lambda$. We exploit the computations from sections 4.2.6 and A.6 in \cite{BP}.   
The optimal trajectory to go from $\eps$ to $a$ is the original ODE, but time reversed, i.e. it follows the ODE 
$\dot{x}_t=\beta_2(x_t)-\beta_1(x_t)$. The running cost is 
$(\beta_2(x_t)-\beta_1(x_t))\log\left(\frac{\beta_2(x_t)}{\beta_1(x_t)}\right)$, 
so the total cost is
\begin{align*}
\int_{T_\eps}^{T_a} \log\left(\frac{\beta_2(x_t)}{\beta_1(x_t)}\right)(\beta_2(x_t)-\beta_1(x_t))dt&=\int_{T_.eps}^{T_a} 
\log\left(\frac{\beta_2(x_t)}{\beta_1(x_t)}\right)\dot{x}_tdt\\
&= \int_\eps^a\log\left(\frac{\beta_2(x)}{\beta_1(x)}\right)dx\\
&\to a+\left(\frac{\gamma}{\lambda}-a\right)\log\left(1-a\frac{\lambda}{\gamma}\right),\ \text{ as }\eps\to0.
\end{align*}
Consequently, from Theorem \ref{LD}, for any $\eta>0$, there exists $t$ and $N$ large enough such that 
\[ \P(Z^N_t\ge a)\le \exp\left\{-N\left[a+\left(\frac{\gamma}{\lambda}-a\right)\log\left(1-a\frac{\lambda}{\gamma}\right)-\eta\right]\right\}.\]
We note that Moderate Deviations resembles much more the Central Limit Theorem than Large Deviations. The fact that the discontinuity in the form of the rate function is exactly at $\alpha=0$ is typical of random variables with light tails. The situation would be quite different with heavy tails, see 
e. g. section VIII.4 in Petrov \cite{P}.

Note however that for small $a$, 
\[ a+\left(\frac{\gamma}{\lambda}-a\right)\log\left(1-a\frac{\lambda}{\gamma}\right)\sim \frac{\lambda a^2}{2\gamma},\]
which is not too surprising, and in a sense reconciles Large Deviations and Moderate Deviations. Were our driving noises Brownian,
then the LD rate function would be quadratic as that of MD, but the LD quasi--potential is the minimal cost when controlling the LLN ODE, while the MD quasi--potential is the minimal cost when controlling the linearized ODE around the endemic equilibrium.

\section{Appendix}
In this Appendix, we establish the following technical result.
\begin{proposition}\label{le:estimexp}
Let $\M$ be a standard Poisson random mesure on $\R^2_+$, and $\MM(dt,du)=\M(dt,du)-dt\,du$ the associated compensated measure. 
If $\varphi$ is an $\R_+$--valued predictable process such that $\int_0^T\varphi_t dt$ has exponential moments of any order, and
$a\in\R$, then for any $0\le t\le T$, 
\[ \E\left[\sup_{0\le s\le t}\exp\left\{a\int_0^s\int_0^{\varphi_r}\MM(dr,du)\right\}\right]\lesssim
\left(\E\exp\left\{(e^{2a}-1-2a)\int_0^t\varphi_sds\right\}\right)^{1/2}\,.\]
\end{proposition}
\bpf 
Consider with $b\ge0$ the process
\begin{equation}\label{eq:X}
X_t=a\int_0^t \int_0^{\varphi_s} \MM(ds,du)-b\int_0^t \varphi_s ds\,.
\end{equation}
It follows from It\^o's formula that
\begin{align*} 
e^{X_t}&=1- b\int_0^te^{X_s}\varphi_sds+a\int_0^t\int_0^{\varphi_s}e^{X_{s-}}\MM(ds,du)\\&\qquad+(e^a-1-a)\int_0^t\int_0^{\varphi_s}e^{X_{s-}}\M(ds,du)\,.
\end{align*}
From Lemma \ref{le:mart} below, $M_t=\int_0^t\int_0^{\varphi_s}e^{X_{s-}}\MM(ds,du)$ is a martingale.
Hence $e^{X}$ is a martingale if $b=(e^a-1-a)$, a submartingale if we replace $=$ by $<$, and a supermartingale if we replace $=$ by $>$. Consequently if $b\ge(e^a-1-a)$, $\E e^{X_t}\le1$.
Now, using first Doob's $L^2$ inequality for submartingales, and later Schwartz's inequality, we have
\begin{align*}
&\E\left[\sup_{0\le s\le t}\exp\left\{a\int_0^s \int_0^{\varphi_r} \MM(dr,du)\right\}\right]\\
&\lesssim\E\exp\left\{a\int_0^t \int_0^{\varphi_s} \MM(ds,du)\right\}\\&=\E\left(\exp\left\{a\int_0^t \int_0^{\varphi_s} \MM(ds,du)-b\int_0^t \varphi_s ds\right\}\exp\left\{b\int_0^t \varphi_s ds\right\}\right)\\
&\le\left(\E\exp\left\{2a\int_0^t \int_0^{\varphi_s} \MM(ds,du)-2b\int_0^t \varphi_s ds\right\}\right)^{1/2}\!\!\!\!\!\!\times\!\!
\left(\E\exp\left\{2b\int_0^t \varphi_s ds\right\}\right)^{1/2}
\end{align*}
If $2b= e^{2a}-1-2a$, it follows from the previous argument that the first factor on the second right hand side is less than or equal to $1$, hence the result follows.
\epf

In order to complete the proof of Proposition \ref{le:estimexp}, we still need to establish
\begin{lemma}\label{le:mart}
The process $\varphi$ satisfying the same assumptions as in Proposition \ref{le:estimexp}, and $X_t$ being given by \eqref{eq:X}, 
$M_t=\int_0^t\int_0^{\varphi_s}e^{X_{s-}}\MM(ds,du)$ is a martingale.
\end{lemma}
\bpf It is plain that $M_t$ is a local martingale, whose predictable quadratic variation is given as
\[ <M>_t=\int_0^te^{2X_s}\varphi_sds\begin{cases}
\le\exp\left\{2a\int_0^t\int_0^{\varphi_s}\M(ds,du)\right\}\int_0^t\varphi_sds,&\text{if $a>0$}\, ;\\
\le\exp\left\{-2(a+b)\int_0^t\varphi_sds\right\}\int_0^t\varphi_sds,&\text{if $a\le0$}\,.
\end{cases}\] 
All we need to show is that the above quantity is integrable. 
It is clearly a consequence of the assumption in case $a<0$. In case $a>0$, the second factor of the right hand side has finite exponential moments, so is square integrable, and all we need to show is that
\begin{equation*}\label{eq:fini}
\E\exp\left\{4a\int_0^t\int_0^{\varphi_s}\M(ds,du)\right\}<\infty.
\end{equation*}
Using It\^o's formula we have \begin{align*}
 Y_t&=\exp\left\{8a\int_0^t\int_0^{\varphi_s}\M(ds,du)-(e^{8a}-1)\int_0^t\varphi_sds\right\}\\&=1+(e^{8a}-1)\int_0^t\int_0^{\varphi_s}Y_{s-}\MM(ds,du).
 \end{align*}
The same computation with $\varphi_s$ replaced by $\varphi^n_s=\varphi_s\wedge n$, and then $Y_s$ replaced by $Y^n_s$
would show that $Y_t^n$ is a martingale satisfying $\E Y^n_t=1$. But $0\le Y^n_t\to Y_t$ a.s., hence Fatou's Lemma implies that $\E Y_t\le 1$. Since
\[ 4a\int_0^t\int_0^{\varphi_s}\M(ds,du)=4a\int_0^t\int_0^{\varphi_s}\M(ds,du)-\frac{e^{8a}-1}{2}\int_0^t\varphi_sds+
\frac{e^{8a}-1}{2}\int_0^t\varphi_sds,\]
it follows from Schwartz's inequality that
\[\E\exp\left\{4a\int_0^t\int_0^{\varphi_s}\M(ds,du)\right\}\le\sqrt{\E Y_t}\sqrt{\E\exp\left\{(e^{8a}-1)\int_0^t\varphi_sds\right\}}\, ,\]
and the result follows from our assumption on $\varphi$.  \epf

\paragraph{Acknowledgement}†It is a pleasure to thank Pierre Petit for an inspiring discussion on moderate deviations.


\begin{thebibliography}{99}
\bibitem{BP}  Tom Britton and Etienne Pardoux, Stochastic Epidemics in a Homogeneous Community, arXiv:1808.0535, submitted.
\bibitem{DZ} Amir Dembo and Ofer Zeitouni, {\em Large Deviations, Techniques and Applications}, 2d ed., Applications of Mathematics {\bf38}, Springer, New York, 1998.
\bibitem{EK} Stewart N. Ethier and Thomas G. Kurtz, {\em Markov processes. Characterization and convergence}, J. Wiley 1986.
\bibitem{FW} Mark~I. Freidlin and Alexander~D. Wentzell.
\newblock {\em Random perturbations of dynamical systems}, 3d ed. Grundlehren des Mathematischen Wissenschaften {\bf260},
\newblock Springer, New York, 2012.
\bibitem{KP} Peter Kratz and Etienne Pardoux, Large deviations for infectious diseases models, in
{\it S\'eminaire de Probabilit\'es XLIX}, C. Donati-Martin, A. Lejay, A. Rouault eds., Lecture Notes in Math. {\bf2215}, pp. 221-327, 2018.
\bibitem{PS} Etienne Pardoux and  Brice Samegni--Kepgnou, Large deviation principle for epidemic models, {\it Journal of Applied Probability}  {\bf54}, 905--920, 2017.  
\bibitem{P} Vasily V. Petrov {\em Sums of Independent Random Variables}, Ergebnisse der Mathematik und ihrer Grenzgebiete {\bf 82}, Springer Verlag, 1975.
\bibitem{PBGM} 
Lev S. Pontryagin, Vladimir G. Boltyanskii,  Revaz V. Gamkrelidze and  Evgenii F. Mishchenko  {\it The mathematical theory of optimal processes}. Transl. by K. N. Trirogoff; ed. by L. W. Neustadt, John Wiley \& Sons, 1962.
\bibitem{SW} Adam Shwartz and Alan Weiss  {\it Large Deviations for Performance Analysis}, Chapman Hall, London, 1995
\bibitem{SW05} Adam Shwartz and Alan Weiss, Large deviations with diminishing rates, {\it Mathematics of Operations Research} {\bf30} 281--310, 2005.

  
\end{thebibliography}
\end{document}